\documentclass[a4paper,reqno,12pt]{amsart} 
\usepackage{newclude}
\usepackage[margin=1.2in]{geometry}
\usepackage{indentfirst}
\usepackage{amsfonts}
\usepackage{amscd}
\usepackage{amsthm}
\usepackage{amssymb}
\usepackage{amsmath}
\usepackage{bm}
\usepackage{thmtools}
\usepackage{mathrsfs}
\usepackage{epsfig}
\usepackage{graphicx}
\usepackage{slashed}
\usepackage[all,cmtip]{xy}
\usepackage{tikz}
\usepackage{tikz-cd}
\usepackage[textsize=footnotesize]{todonotes}
\usepackage{extpfeil}
\usetikzlibrary{decorations.markings,shapes.geometric}
\usepackage{enumitem}
\usepackage{setspace}
\usepackage{upgreek}
\usepackage{verbatim}
\usepackage[bookmarksdepth=2,linktoc=page,colorlinks,linkcolor={red!80!black},citecolor={red!80!black},urlcolor={blue!80!black}]{hyperref}

\usetikzlibrary{matrix,arrows,positioning}

\def\itemNum$#1${\item $\displaystyle#1$
   \hfill\refstepcounter{equation}(\theequation)}

\makeatletter
\providecommand\@dotsep{5}
\renewcommand{\listoftodos}[1][\@todonotes@todolistname]{%
  \@starttoc{tdo}{#1}}
\makeatother

\newtheorem{Lem}{Lemma}[section]
\newtheorem{Prop}[Lem]{Proposition}

\newtheorem*{Def}{Definition}

\theoremstyle{plain}
\newtheorem{Thm}[Lem]{Theorem}

\newtheorem{Cor}[Lem]{Corollary}

\newtheorem{thm}{Theorem}

\theoremstyle{definition}
\declaretheorem[numbered=no,name=Example,qed={\lower-0.3ex\hbox{$\triangleleft$}}]{Ex}
\newtheorem{Rem}[Lem]{Remark}

\def\angleser#1{\langle\!\langle#1\rangle\!\rangle} 

\DeclareMathOperator{\Hom}{Hom}
\DeclareMathOperator{\HHom}{\mathcal{H}om}
\DeclareMathOperator{\End}{End}

\DeclareMathOperator{\Aut}{Aut}
\DeclareMathOperator{\Sym}{Sym}
\DeclareMathOperator{\Spec}{Spec}
\DeclareMathOperator{\Proj}{Proj}

\DeclareMathOperator{\Pic}{Pic}
\DeclareMathOperator{\SPic}{Pic^{\text{\textnormal{sym}}}}

\DeclareMathOperator{\Vect}{\mathsf{Vect}}

\DeclareMathOperator{\Modu}{\mhyphen\mathsf{Mod}}
\DeclareMathOperator{\proj}{\mhyphen\mathsf{proj}}
\DeclareMathOperator{\modu}{\mhyphen\mathsf{mod}}

\DeclareMathOperator{\LP}{\mathsf{LP}}
\DeclareMathOperator{\LF}{\mathsf{LF}}
\DeclareMathOperator{\Frac}{Frac}
\DeclareMathOperator{\lf}{lf}
\DeclareMathOperator{\vect}{vect}

\newcommand{\im}{\textup{im}}
\newcommand{\coker}{\text{\textnormal{coker}}}
\mathchardef\mhyphen="2D
\newcommand{\git}{/\!\!/}

\newcommand{\op}{\text{\textnormal{op}}}
\newcommand{\id}{\text{\textnormal{id}}}

\newcommand{\nor}{\text{\textnormal{n}}}

\newcommand{\rpc}{\text{\textnormal{rpc}}}
\newcommand{\rc}{\text{\textnormal{rc}}}

\newcommand{\A}{\mathbb{A}}
\renewcommand{\P}{\mathbb{P}}
\newcommand{\Z}{\mathbb{Z}}
\newcommand{\cE}{\mathcal{E}}
\newcommand{\cF}{\mathcal{F}}
\newcommand{\cL}{\mathcal{L}}

\newcommand{\cO}{\mathcal{O}}

\newcommand{\fp}{\mathfrak{p}}
\newcommand{\fm}{\mathfrak{m}}
\newcommand{\Fun}{{\mathbb{F}_1}}
\newcommand{\gen}[1]{\langle#1\rangle}
\newcommand{\res}{\textup{res}}
\newcommand{\pr}{\textup{pr}}

\makeatletter
\newbox\xrat@below
\newbox\xrat@above
\newcommand{\xrightarrowtail}[2][]{%
  \setbox\xrat@below=\hbox{\ensuremath{\scriptstyle #1}}%
  \setbox\xrat@above=\hbox{\ensuremath{\scriptstyle #2}}%
  \pgfmathsetlengthmacro{\xrat@len}{max(\wd\xrat@below,\wd\xrat@above)+.6em}%
  \mathrel{\tikz [>->,baseline=-.75ex]
                 \draw (0,0) -- node[below=-2pt] {\box\xrat@below}
                                node[above=-2pt] {\box\xrat@above}
                       (1.8*\xrat@len,0) ;}}
\renewcommand{\xtwoheadrightarrow}[2][]{%
  \setbox\xrat@below=\hbox{\ensuremath{\scriptstyle #1}}%
  \setbox\xrat@above=\hbox{\ensuremath{\scriptstyle #2}}%
  \pgfmathsetlengthmacro{\xrat@len}{max(\wd\xrat@below,\wd\xrat@above)+.6em}%
  \mathrel{\tikz [->>,baseline=-.75ex]
                 \draw (0,0) -- node[below=-2pt] {\box\xrat@below}
                                node[above=-2pt] {\box\xrat@above}
                       (1.8*\xrat@len,0) ;}}
\makeatother

\title[$K(X_{\slash \mathbb{F}_1})$ and $GW(X_{\slash \mathbb{F}_1})$]{Algebraic $K$-theory and Grothendieck--Witt theory of monoid schemes}

\author[J.\,N. Eberhardt]{Jens Niklas Eberhardt}
\address{Max Planck Institute for Mathematics\\
Vivatsgasse 7\\
53111 Bonn, Germany}
\email{mail@jenseberhardt.com}

\author[O. Lorscheid]{Oliver Lorscheid}
\address{Instituto Nacional de Matem\'{a}tica Pura e Aplicada, Rio de Janeiro, Brazil}
\email{oliver@impa.br}

\author[M.\,B. Young]{Matthew B. Young}
\address{Department of Mathematics and Statistics \\ Utah State University\\
Logan, Utah 84322 \\ USA}
\email{matthew.young@usu.edu}

\date{\today}


\begin{document}

\begin{abstract}
We study the algebraic $K$-theory and Grothendieck--Witt theory of proto-exact categories of vector bundles over monoid schemes. Our main results are the complete description of the algebraic $K$-theory space of an integral monoid scheme $X$ in terms of its Picard group $\Pic(X)$ and pointed monoid of regular functions $\Gamma(X, \mathcal{O}_X)$ and a description of the Grothendieck--Witt space of $X$ in terms of an additional involution on $\Pic(X)$. We also prove space-level projective bundle formulae in both settings.
\end{abstract}

\maketitle

\begin{small} \tableofcontents \end{small}

\setcounter{footnote}{0}

\section*{Introduction}
\addtocontents{toc}{\protect\setcounter{tocdepth}{1}}

In this paper we study the algebraic $K$-theory and Grothendieck--Witt theory of monoid schemes. Our main results state that, under mild assumptions on the scheme $X$, these spaces are determined by simple algebraic invariants of $X$. For the $K$-theory space, these invariants are the Picard group $\Pic(X)$ and the group of invertible regular functions $\Gamma(X, \mathcal{O}_X)^{\times}$. For the Grothendieck--Witt space, the additional data of a set-theoretic involution
of $\Pic(X)$ is required.

Monoid schemes form the core of algebraic geometry over the elusive field $\Fun$ with one element \cite{tits1957}, \cite{soule2004}, in the sense that every other approach to $\Fun$-schemes contains monoid schemes as a full subcategory. In the other direction, monoid schemes can be seen as a direct generalization of toric geometry and Kato fans of logarithmic schemes; see \cite{kato1994}, \cite{deitmar2005}, \cite{connes2010}, \cite{chu2012}, \cite{cortinas2015} among others. The central position of monoid schemes within $\Fun$-geometry is confirmed by the multiple links to other disciplines, such as Weyl groups as algebraic groups over $\Fun$ \cite{lorscheid2012}, computational methods for toric geometry \cite{cortinas2014}, \cite{cortinas2015}, \cite{flores2017}, a framework for tropical scheme theory \cite{giansiracusa2016}, applications to representation theory \cite{szczesny2018}, \cite{jun2020} and, last but not least, stable homotopy theory as $K$-theory over $\Fun$ \cite{deitmar2006}, \cite{chu2012}, the theme on which we dwell in this paper.

For the purpose of this introduction, we provide the reader with the following suggestive description: a monoid scheme is a topological space together with a sheaf of commutative pointed monoids which is locally isomorphic to the spectrum of a commutative pointed monoid. Here, a pointed monoid is a monoid $A$ with an absorbing element $0$, that is, an element satisfying $0\cdot a=0$ for all $a \in A$.

Just as algebraic $K$-theory is an algebraic analogue of complex topological $K$-theory, Grothendieck--Witt theory is an algebraic analogue of Atiyah's topological $KR$-theory \cite{atiyah1966}. A key feature of $KR$-theory is that it generalizes complex, real and quaternionic topological $K$-theory. Grothendieck--Witt theory enjoys a similar status in the algebraic setting. For example, whereas the algebraic $K$-theory of a scheme $X$ studies algebraic vector bundles on $X$, the Grothendieck--Witt theory of $X$ studies algebraic vector bundles with non-degenerate bilinear form on $X$ or, equivalently, orthogonal or symplectic vector bundles on $X$, depending on the symmetry of the pairing. Grothendieck--Witt theory plays a fundamental role in Karoubi's formulation and proof of topological and algebraic Bott periodicity and study of the homology of orthogonal and symplectic groups \cite{karoubi1973}, \cite{karoubi1980}, \cite{karoubi1980b}. Recently, much effort has been devoted to developing the Grothendieck--Witt theory of schemes; see, for example, \cite{fasel2009}, \cite{schlichting2010}, \cite{schlichting2015}, \cite{karoubi2016}, \cite{karoubi2017}, \cite{karoubi2020}. This paper provides the first results in the development of these ideas for monoid schemes.

In Section \ref{sec:prelim} we recall relevant categorical and $K$-theoretic background. We work in the setting of proto-exact categories, a non-additive generalization of Quillen's exact categories introduced by Dyckerhoff and Kapranov \cite{dyckerhoff2019}. This is a convenient setting for both $K$-theory and Grothendieck--Witt theory and was developed by the authors in \cite{eberhardtLorscheidYoung2020a}.
The first main result of \cite{eberhardtLorscheidYoung2020a} is a Group Completion Theorem for the $K$-theory of uniquely split proto-exact categories. The second main result of \cite{eberhardtLorscheidYoung2020a} is a description of the Grothendieck--Witt space $\mathcal{GW}^Q(\mathcal{A})$ of a uniquely split proto-exact category with duality $\mathcal{A}$ satisfying additional mild assumptions, defined using the hermitian $Q$-construction, in terms of the group completion of the groupoid of hyperbolic forms and the monoidal groupoid of isotropically simple symmetric forms in $\mathcal{A}$. The second result can be seen as playing the role of the Group Completion Theorem in the Grothendieck--Witt theory of uniquely split proto-exact categories with duality.

In Section \ref{sec:KGWMonoids} we study the $K$-theory and Grothendieck--Witt theory of proto-exact categories of modules over pointed monoids or, geometrically, vector bundles over non-commutative affine monoid schemes. Let $A$ be a (not necessarily commutative) pointed monoid. The category $A \proj$ of finitely generated projective left $A$-modules has a uniquely split proto-exact structure (Lemma \ref{lem:AProjSplit}), and so fits into the framework of \cite{eberhardtLorscheidYoung2020a}. The category $A \proj$ is particularly simple when $A$ is integral,\footnote{More generally, $A$ need only be right partially cancellative. We work at this level of generality in the body of the paper.}
in which case all projective $A$-modules are free. This, together with the Group Completion Theorem, allows for an explicit description of the $K$-theory space $\mathcal{K}(A \proj)$ 
in terms of the group of units $A^{\times}$, thereby giving a `$Q=+$' theorem in this setting. In this way, we extend earlier results of Deitmar \cite{deitmar2006} and establish some unproven claims of Chu--Morava \cite{chu2010}. Our result is as follows.
\begin{thm}[{Theorem \ref{thm:KThyRPC}}]
\label{thm:KThyRPCIntro}
Let $A$ be an integral pointed monoid. Then there is a homotopy equivalence
\[
\mathcal{K}(A \proj)
\simeq
\mathbb{Z} \times B (A^{\times} \wr \Sigma_{\infty})^+,
\]
where $A^{\times} \wr \Sigma_{\infty}$ is the infinite wreath product $\displaystyle \varinjlim_n\, ((A^{\times})^n \rtimes \Sigma_n)$.
\end{thm}

Next, we study the Grothendieck--Witt theory of projective $A$-modules, where our results are new. Unlike the case of rings, the category $A \proj$ does not admit an exact duality structure. In particular, the functor $\Hom_{A \proj}(-,A)$ is poorly behaved. For this reason, we restrict attention to the non-full proto-exact subcategory $A \proj^{\nor} \subset A \proj$ of normal morphisms, that is, $A$-module homomorphisms whose non-empty fibres over non-basepoints are singletons. We prove in Lemma \ref{lem:dualFreeAMod} that, if $A$ is integral, then $A \proj^{\nor}$ admits a duality structure. From this point of view, normal morphisms are essential to Grothendieck--Witt theory. We remark that the $K$-theory of $A \proj^{\nor}$ and $A \proj$ coincides. For earlier appearances of normal morphisms in $\Fun$-geometry, see \cite{chu2012}, \cite{szczesny2014}, \cite{mbyoung2021}.
The main results of this section determine the Grothendieck--Witt spaces $\mathcal{GW}^{\oplus}(A \proj^{\nor})$ and $\mathcal{GW}^Q(A \proj^{\nor})$. In particular, the former space can be described in terms of ``infinite orthogonal/symplectic groups over $\Fun$". In this way, we obtain an $\Fun$-analogue of Karoubi's results on the hermitian $K$-theory of rings \cite{karoubi1973}. A simplified version of this result is as follows.

\begin{thm}[{Theorem \ref{thm:GWThyMonoid}}]
\label{thm:GWThyMonoidIntro}
Let $A$ be an integral pointed monoid with $A^{\times} = \{1\}$. Then there is a homotopy equivalence
\[
\mathcal{GW}^{\oplus}(A \proj^{\nor})
\simeq
\mathbb{Z}^2 \times B(\mathbb{Z}/2 \wr \Sigma_{\infty})^+.
\]
\end{thm}

By applying the results of \cite{eberhardtLorscheidYoung2020a}, we can use Theorem \ref{thm:GWThyMonoidIntro} to describe the weak homotopy type of $\mathcal{GW}^Q(A \proj^{\nor})$.

Having treated the local theory, we turn in Sections \ref{sec:KThySch} and \ref{sec:GWThySch} to the global theory of monoid schemes. Following earlier approaches \cite{huttemann2002}, \cite{deitmar2006}, a general definition of the $K$-theory of a monoid scheme $X$ was given in \cite{chu2012}, where a proto-exact category of vector bundles $\Vect(X)$ and their normal $\mathcal{O}_X$-module homomorphisms was defined. We point out that this is not the only approach to the $K$-theory of monoid schemes; see \cite{haesemeyer2019} for a recent alternative. In Section \ref{sec:KThySch} we bring the approach of \cite{chu2012} to its natural conclusion by explicitly describing the $K$-theory space $\mathcal{K}(X) := \mathcal{K}(\Vect(X))$. The key structural result is Proposition \ref{prop:VectDecomp}, which exhibits an extremely simple non-full proto-exact subcategory $\angleser{\mathcal{O}_X}[\Pic(X)]$
of $\Vect(X)$ whose $K$-theory space is homotopy equivalent to $\mathcal{K}(X)$. Here, $\angleser{\mathcal{O}_X}$ is the category whose objects are isomorphic to $\mathcal{O}_X^{\oplus n}$, $n \in \mathbb{Z}_{\geq 0}$, together with all normal $\mathcal{O}_X$-module homomorphisms between them. The category $\angleser{\mathcal{O}_X}[\Pic(X)]$ can then be seen as the group algebra of $\Pic(X)$ with coefficients in $\angleser{\mathcal{O}_X}$; see Section \ref{sec:LPsheaves} for a precise definition. Proposition \ref{prop:VectDecomp} fails for the exact category of vector bundles over a field and is the source of the relative strength of the following result.

\begin{thm}[{Theorem \ref{thm:KThySchemes}}]
\label{thm:KThySchemesIntro}
Let $X$ be an integral monoid scheme. Then there is a homotopy equivalence
\[
\mathcal{K}(X)
\simeq
\sideset{}{'}\prod_{\mathcal{M} \in \Pic(X)} \mathbb{Z} \times B (\Gamma(X, \mathcal{O}_X)^{\times} \wr \Sigma_{\infty})^+,
\]
where $\prod^{\prime}$ is the restricted product of pointed topological spaces.
\end{thm}

Turning to Grothendieck--Witt theory, let $\mathcal{L}$ be a line bundle on an integral monoid scheme $X$. The integrality assumption on $X$ ensures the existence of a duality structure $(P^{\mathcal{L}}, \Theta^{\mathcal{L}})$ on $\Vect(X)$. Write $\mathcal{GW}(X;\mathcal{L}) := \mathcal{GW}(\Vect(X), P^{\mathcal{L}}, \Theta^{\mathcal{L}})$ for the associated Grothendieck--Witt space, defined either via the hermitian $Q$-construction or group completion. The duality structure $(P^{\mathcal{L}}, \Theta^{\mathcal{L}})$ is compatible with the subcategory $\angleser{\mathcal{O}_X}[\Pic(X)]$, which leads to a complete description of $\mathcal{GW}(X;\mathcal{L})$ in terms of $\Pic(X)$, together with its set-theoretic $\mathbb{Z}/2$-action determined by $\mathcal{L}$, and the pointed monoid $\Gamma(X, \mathcal{O}_X)$.

\begin{thm}[{Theorem \ref{thm:GWThySchemes}}]
\label{thm:GWThySchemesIntro}
Let $\mathcal{L}$ be a line bundle on an integral monoid scheme $X$. Then there is a natural homotopy equivalence
\[
\mathcal{GW}(X;\mathcal{L})
\simeq
\sideset{}{'}\prod_{\mathcal{M} \in \Pic(X)^{P^{\mathcal{L}}}} \mathcal{GW}(\Gamma(X, \mathcal{O}_X) \proj^{\nor})
\times
\sideset{}{'}\prod_{\mathcal{M} \in \Pic(X)^* \slash P^{\mathcal{L}}} \mathcal{K}(\Gamma(X, \mathcal{O}_X) \proj^{\nor}),
\]
where $\Pic(X)^{P^{\mathcal{L}}}$ denotes the fixed point set of $\Pic(X)$ under the $\mathbb{Z} \slash 2$-action determined by $\mathcal{L}$ and $\Pic(X)^* \slash P^{\mathcal{L}}$ is the quotient of the complement $\Pic(X) \setminus \Pic(X)^{P^{\mathcal{L}}}$.
\end{thm}

In particular, this result, together with Theorems \ref{thm:KThyRPCIntro} and \ref{thm:GWThyMonoidIntro}, leads to an explicit description of $\mathcal{GW}^{\oplus}(X;\mathcal{L})$. The space $\mathcal{GW}^Q(X;\mathcal{L})$ can then be described using the results of \cite{eberhardtLorscheidYoung2020a}.

As an application of our results, we prove space-level projective bundle formulae for $K$-theory and Grothendieck--Witt theory, giving analogues of well-known results for schemes over fields \cite{quillen1973}, \cite{walter2003}, \cite{schlichting2017}, \cite{rohrbach2020}. In our setting, the key background results are Theorem \ref{thm:ProjLineBundles} and Lemma \ref{lem:equivarOfPhi}, which give a ($\mathbb{Z}/2$-equivariant) description of the Picard group of a projective bundle in terms of that of the base. Notably, our proof relies on different arguments than the classical proof, since sheaf cohomology is not available in the $\Fun$-setting. Instead, we use particular properties of monoids schemes whose analogues for schemes over fields fail to hold. 
The projective bundle formula for Grothendieck--Witt theory is as follows; for $K$-theory, see Theorem \ref{thm:projBunSpaces}.

\begin{thm}[{Theorem \ref{thm:projBundGWSpaces}}]
\label{thm:projBundGWSpacesIntro}
Let $\mathcal{E}$ be a vector bundle on an integral monoid scheme $X$ with associated projective bundle $\pi: \mathbb{P} \mathcal{E} \rightarrow X$ and $\cL$ a line bundle on $X$. Then there is a homotopy equivalence
\[
\mathcal{GW}(\mathbb{P}\mathcal{E}; \pi^*\mathcal{L})
\simeq
\mathcal{GW}(X; \mathcal{L}) \times 
\sideset{}{'}\prod_{(\mathcal{M},i) \in (\Pic(X) \times \mathbb{Z}^*) \slash \langle (P^{\mathcal{L}},-1) \rangle} \mathcal{K}(\Gamma(X, \mathcal{O}_X) \proj^{\nor}).
\]
\end{thm}

\subsubsection*{Acknowledgements}
The authors thank Marco Schlichting for helpful correspondence. All three authors thank the Max Planck Institute for Mathematics in Bonn for its hospitality and financial support.


\section{Background material}
\label{sec:prelim}

In this section, we recall necessary background material on proto-exact categories and their $K$-theory and Grothendieck--Witt theory.

\subsection{Proto-exact categories}
\label{sec:cat}

Let $\mathcal{A}$ be a proto-exact category, as defined in \cite[\S 2.4]{dyckerhoff2019}. In particular, the category $\mathcal{A}$ has a zero object $0 \in \mathcal{A}$ and two distinguished classes of morphisms, called inflations (or admissible monics) and deflations (or admissible epics) and denoted by $\rightarrowtail$ and $\twoheadrightarrow$, respectively. An admissible square in $\mathcal{A}$ is a bicartesian square of the form
\[
\begin{tikzpicture}[baseline= (a).base]
\node[scale=1] (a) at (0,0){
\begin{tikzcd}
U \arrow[two heads]{d} \arrow[tail]{r} & V \arrow[two heads]{d}\\
W \arrow[tail]{r} & X
\end{tikzcd}
};
\end{tikzpicture}
\]
Conflations (or admissible short exact sequences) in $\mathcal{A}$ are admissible squares as above with $W=0$ which, for ease of notation, we denote by $U \rightarrowtail V \twoheadrightarrow X$.

A functor between proto-exact categories is called proto-exact if it sends admissible squares to admissible squares. In particular, proto-exact functors send conflations to conflations.

The proto-exact categories of interest in this paper have a weak analogue of an additive structure, axiomatized as follows.

\begin{Def}[{\cite[\S 1.1]{eberhardtLorscheidYoung2020a}}]
An exact direct sum on a proto-exact category $\mathcal{A}$ is a symmetric monoidal structure $\oplus$ on $\mathcal{A}$ such that $0$ is the monoidal unit and $\oplus$ is a proto-exact functor. Moreover, the following additional axioms are required to hold, where we set $i_U: U \xrightarrowtail[]{\id_U \oplus 0_{0 \rightarrowtail V}} U \oplus V$ and $\pi_U: U \oplus V \xtwoheadrightarrow[]{\id_U \oplus 0_{V \twoheadrightarrow 0}} U$ for objects $U,V \in \mathcal{A}$.
\begin{enumerate}[wide,labelwidth=!, labelindent=0pt,label=(\roman*)]
\item \label{ax:restrBij} The map
\[
	\Hom_{\mathcal{A}}(U \oplus V, W) \rightarrow \Hom_{\mathcal{A}}(U,W) \times \Hom_{\mathcal{A}}(V,W),
\qquad
f \mapsto (f \circ i_U, f \circ i_V)
\]
is an injection for all $U,V,W \in \mathcal{A}$.

\item \label{ax:sectSplit} Let $U \xrightarrowtail[]{i} X \xtwoheadrightarrow[]{\pi} V$ be a conflation. For each section $s$ of $\pi$, there exists a unique isomorphism $\phi$ which makes the following diagram commute:
\[
\begin{tikzpicture}[baseline= (a).base]
\node[scale=1] (a) at (0,0){
\begin{tikzcd}
& X & \\
U \arrow[tail]{r}[below]{i_U} \arrow[tail]{ur}[above]{i}& U \oplus V \arrow{u}[left]{\phi}& \arrow[tail]{l}[below]{i_V} \arrow{ul}[above]{s} V.
\end{tikzcd}
};
\end{tikzpicture}
\]
\end{enumerate}
Moreover, the obvious axioms dual to \ref{ax:restrBij} and \ref{ax:sectSplit}, with the maps $\pi_{(-)}$ appearing in place of $i_{(-)}$, are required to hold.
\end{Def}

A functor between proto-exact categories with exact direct sum is called exact if it is proto-exact and $\oplus$-monoidal.

Let $\mathcal{A}$ be a proto-exact category with exact direct sum. A commutative diagram
\[
\begin{tikzpicture}[baseline= (a).base]
\node[scale=1] (a) at (0,0){
\begin{tikzcd}
U \arrow[tail]{r}[above]{i} \arrow[d,equal] & X \arrow[two heads]{r}[above]{\pi} & V \arrow[d,equal]\\
U \arrow[tail]{r}[below]{i_U} & U \oplus V \arrow{u}[left]{\phi} \arrow[two heads]{r}[below]{\pi_V}& V
\end{tikzcd}
};
\end{tikzpicture}
\]
with $\phi$ an isomorphism is called a splitting of the conflation $U \xrightarrowtail[]{i} X \xtwoheadrightarrow[]{\pi} V$.

\begin{Def}[{\cite[\S 1.1]{eberhardtLorscheidYoung2020a}}]
A proto-exact category with exact direct sum is called
\begin{enumerate}[wide,labelwidth=!, labelindent=0pt,label=(\roman*)]
\item uniquely split if every conflation admits a unique splitting, and

\item combinatorial if, for each inflation $i: U \rightarrowtail X_1 \oplus X_2$, there exist inflations $i_k: U_k \rightarrowtail X_k$, $k=1,2$, and an isomorphism $f: U \rightarrow U_1 \oplus U_2$ such that $i = (i_1 \oplus i_2) \circ f$. Moreover, the obvious dual axiom involving maps $\pi_k$, $k=1,2$, is required to hold.
\end{enumerate}
\end{Def}

\subsection{Algebraic $K$-theory of proto-exact categories}
\label{sec:algKProtoEx}

Let $\mathcal{A}$ be a proto-exact category. The $Q$-construction of $\mathcal{A}$ can be defined as for exact categories \cite[\S 2]{quillen1973}, yielding a category $Q(\mathcal{A})$. See also \cite[\S 2.1]{eberhardtLorscheidYoung2020a}. The $K$-theory space of $\mathcal{A}$ is then $\mathcal{K}(\mathcal{A}) = \Omega B Q(\mathcal{A})$, where $BQ(\mathcal{A})$ is pointed by $0 \in Q(\mathcal{A})$, and the $K$-theory groups are
\[
K_i(\mathcal{A}) = \pi_i \mathcal{K}(\mathcal{A}),
\qquad
i \geq 0.
\]

\begin{Lem}
\label{lem:reflectExactIso}
Let $\mathcal{A}$ and $\mathcal{B}$ be proto-exact categories and $F: \mathcal{A} \rightarrow \mathcal{B}$ an essentially surjective proto-exact functor which is bijective on inflations and deflations. Then the induced map $\mathcal{K}(F) : \mathcal{K}(\mathcal{A}) \rightarrow \mathcal{K}(\mathcal{B})$ is a homotopy equivalence.
\end{Lem}

\begin{proof}
To begin, note that $F$ is conservative. Indeed, a morphism in a proto-exact category is an isomorphism if and only if it is an inflation and a deflation.

Since $F$ is proto-exact, there is an induced functor $Q(F): Q(\mathcal{A}) \rightarrow Q(\mathcal{B})$. Essential surjectivity of $F$ implies that of $Q(F)$. Moreover, $Q(F)$ is full (resp. faithful) because $F$ is surjective on inflations and deflations (resp. conservative and injective on inflations and deflations). Hence, $Q(F)$ is an equivalence and the associated map $\mathcal{K}(F)$ is a homotopy equivalence.
\end{proof}

Let now $(\mathcal{A}, \oplus)$ be a symmetric monoidal category. The maximal groupoid $\mathcal{S} \subset \mathcal{A}$ inherits a symmetric monoidal structure. Following \cite[Page 222]{grayson1976}, the direct sum $K$-theory space of $\mathcal{A}$ is the group completion of $B\mathcal{S}$:
\[
\mathcal{K}^{\oplus}(\mathcal{A})
=
B (\mathcal{S}^{-1} \mathcal{S}).
\]

We have the following proto-exact analogue of Quillen's Group Completion Theorem \cite{grayson1976}.

\begin{Thm}[{\cite[Theorem 2.2]{eberhardtLorscheidYoung2020a}}]
\label{thm:KComparison}
Let $\mathcal{A}$ be a uniquely split proto-exact category. Then there is a homotopy equivalence $\mathcal{K}(\mathcal{A}) \simeq \mathcal{K}^{\oplus}(\mathcal{A})$.
\end{Thm}

\begin{Rem}
The construction of the space $\mathcal{K}(\mathcal{A})$ can be refined to produce a  connective spectrum $\mathbf{K}(\mathcal{A})$; see \cite[Remark IV.6.5.1, \S IV.8.5.5]{weibel2013}. While $\mathcal{K}(\mathcal{A})$ and $\mathbf{K}(\mathcal{A})$ have the same homotopy groups, the space $\mathbf{K}(\mathcal{A})$ has many technical advantages. For example, a functor $\otimes: \mathcal{A}\times \mathcal{A}\rightarrow \mathcal{A}$ which is biexact in the sense of \cite[Definition IV.6.6]{weibel2013} induces a pairing of spectra $\mathbf{K}(\mathcal{A})\wedge\mathbf{K}(\mathcal{A})\rightarrow \mathbf{K}(\mathcal{A})$. This gives $K_{\bullet}(\mathcal{A}) = \bigoplus_{i \geq 0} K_i(\mathcal{A})$ the structure of commutative $\mathbb{Z}_{\geq 0}$-graded ring if $\otimes$ is symmetric monoidal.
\end{Rem}

\subsection{Proto-exact categories with duality}
\label{sec:catWD}

For a detailed introduction to proto-exact categories with duality, the reader is referred to \cite[\S 2]{schlichting2010}, \cite[\S 1.2]{eberhardtLorscheidYoung2020a}.

A category with duality is a triple $(\mathcal{A},P, \Theta)$ (often simply $\mathcal{A}$) consisting of a category $\mathcal{A}$, a functor $P: \mathcal{A}^{\op} \rightarrow \mathcal{A}$ and a natural isomorphism $\Theta: \id_{\mathcal{A}} \Rightarrow P \circ P^{\op}$ which satisfies
\begin{equation}
\label{eq:natTransConstr}
P(\Theta_U) \circ \Theta_{P(U)} = \id_{P(U)},
\qquad
U \in \mathcal{A}.
\end{equation}
If $\mathcal{A}$ is proto-exact and $P$ is proto-exact, then $\mathcal{A}$ is a proto-exact category with duality. We henceforth restrict attention to this case.

A symmetric form in $\mathcal{A}$ is an isomorphism $\psi_M: M \rightarrow P(M)$ which satisfies $P(\psi_M) \circ \Theta_M = \psi_M$. An isometry $\phi: (M,\psi_M) \rightarrow (N, \psi_N)$ is an isomorphism $\phi: M \rightarrow N$ which satisfies $\psi_M = P(\phi) \circ \psi_N \circ \phi$. The groupoid of symmetric forms and their isometries is $\mathcal{A}_h$.

Let $(M,\psi_M)$  be a symmetric form. An inflation $i: U \rightarrowtail M$ is called isotropic if $P(i) \circ \psi_M \circ i$ is zero and $U \rightarrow U^{\perp} := \ker(P(i) \circ \psi_M)$ is an inflation. In this case, the reduction $M \git U: = U^{\perp} \slash U$ inherits a symmetric morphism $\psi_{M \git U} : M \git U \rightarrow P(M \git U)$, which we assume to be an isomorphism; this is the Reduction Assumption of \cite[\S 3.4]{mbyoung2018b}.  A symmetric form $(M, \psi_M)$ is called metabolic if it has a Lagrangian, that is, an isotropic subobject $U \rightarrowtail M$ with $U = U^{\perp}$,
and is called isotropically simple if it has no non-zero isotropic subobjects.

If $\mathcal{A}$ has an exact direct sum, then we require that $P$ be exact and $\Theta$ be $\oplus$-monoidal. In this case, $\mathcal{A}_h$ is a symmetric monoidal groupoid. Given an object $U \in \mathcal{A}$, the pair
\[
\left(
H(U) = U \oplus P(U), \psi_{H(U)} = \begin{psmallmatrix} 0 & \id_{P(U)} \\ \Theta_U & 0 \end{psmallmatrix}
\right)
\]
is a symmetric form in $\mathcal{A}$, called the hyperbolic form on $U$. The assignment $U \mapsto H(U)$ extends to a functor $H: \mathcal{S} \rightarrow \mathcal{A}_h$ where $\mathcal{S}$ is the maximal grupoid in $\mathcal{A}$.
A symmetric form which is isometric to $(H(U), \psi_{H(U)})$ for some $U \in \mathcal{A}$ is called hyperbolic.

\begin{Lem}[{\cite[Lemma 1.7]{eberhardtLorscheidYoung2020a}}]
\label{lem:noMeta}
A metabolic form in a uniquely split proto-exact category with duality is hyperbolic.
\end{Lem}

\begin{Ex}
Let $\mathcal{A}$ be a category. The triple $(H(\mathcal{A}), P, \id_{\id_{H(\mathcal{A})}})$, where $H(\mathcal{A}) = \mathcal{A} \times \mathcal{A}^{\op}$ and $P(U,V) = (V,U)$, is called the hyperbolic category with duality on $\mathcal{A}$. If $\mathcal{A}$ is proto-exact, then so too is $H(\mathcal{A})$ and an exact direct sum on $\mathcal{A}$ induces one on $H(\mathcal{A})$.
\end{Ex}

A form functor $(T,\eta): (\mathcal{A}, P, \Theta) \rightarrow (\mathcal{B}, Q, \Xi)$ between categories with duality is a functor $T: \mathcal{A} \rightarrow \mathcal{B}$ and a natural transformation $\eta: T \circ P \Rightarrow Q \circ T^{\op}$ which makes the diagram\footnote{For legibility, we have written $P^2$ in place of $P \circ P^{\op}$, and so on.}
\[
\begin{tikzpicture}[baseline= (a).base]
\node[scale=1.0] (a) at (0,0){
\begin{tikzcd}[column sep=small]
T(U) \arrow{r}[above]{\Xi_{T(U)}} \arrow{d}[left]{T(\Theta_U)} & [2.5em] Q^2T(U) \arrow{d}[right]{Q(\eta_U)}  \\
[1em] TP^2(U) \arrow{r}[below]{\eta_{P(U)}} & QTP(U)
\end{tikzcd}
};
\end{tikzpicture}
\]
commute for each $U \in \mathcal{A}$. The form functor is called non-singular if $\eta$ is a natural isomorphism and is called an equivalence if, moreover, $T$ is an equivalence.

\subsection{Grothendieck--Witt theory of proto-exact categories}
\label{sec:GWProtoEx}

Let $\mathcal{A}$ be a proto-exact category with duality. The hermitian $Q$-construction of $\mathcal{A}$ can be defined as for exact categories with duality \cite[\S 4.1]{schlichting2010}, yielding a category $Q_h(\mathcal{A})$. See also \cite[\S 3.1]{eberhardtLorscheidYoung2020a}. Forgetting symmetric forms defines a functor $F: Q_h(\mathcal{A}) \rightarrow Q(\mathcal{A})$. The Grothendieck--Witt space $\mathcal{GW}^Q(\mathcal{A}$) is the homotopy fibre of $BF: B Q_h(\mathcal{A}) \rightarrow B Q(\mathcal{A})$ over $0$ and the Grothendieck--Witt groups are
\[
GW^Q_i(\mathcal{A}) = \pi_i \mathcal{GW}^Q(\mathcal{A}),
\qquad
i \geq 0.
\]
Despite the name, without further assumptions, $GW^Q_0(\mathcal{A})$ is in fact only a pointed set.  If, however, $\mathcal{A}$ has an exact direct sum, as will always be the case, then $GW^Q_0(\mathcal{A})$ is a commutative monoid. The Witt groups are defined by
\[
W^Q_i(\mathcal{A}) = \coker \left( K_i(\mathcal{A}) \xrightarrow[]{H_*} GW^Q_i(\mathcal{A}) \right),
\qquad
i \geq 0,
\]
where $H_*$ is induced by the map $\mathcal{K}(A) \rightarrow \mathcal{GW}^Q(\mathcal{A})$. As for $GW^Q_0(\mathcal{A})$, in general $W^Q_0(\mathcal{A})$ is only a commutative monoid.

\begin{Prop}
\label{prop:GWHyper}
Let $\mathcal{A}$ be a proto-exact category with associated hyperbolic category $H (\mathcal{A})$. Then there is a homotopy equivalence $\mathcal{GW}^Q(H(\mathcal{A})) \xrightarrow[]{\sim} \mathcal{K}(\mathcal{A})$.
\end{Prop}

\begin{proof}
The proof of the corresponding result in the exact setting \cite[Proposition 4.7]{schlichting2010} carries over.
\end{proof}

\begin{Prop}
\label{prop:GWFunctoriality}
A non-singular 
proto-exact form functor $(T,\eta): (\mathcal{A}, P, \Theta) \rightarrow (\mathcal{B}, Q, \Xi)$ induces a continuous map $\mathcal{GW}^Q(T, \eta): \mathcal{GW}^Q(\mathcal{A}) \rightarrow \mathcal{GW}^Q(\mathcal{B})$. Moreover, if $(T,\eta)$ and $(T^{\prime}, \eta^{\prime})$ are naturally isomorphic, then $\mathcal{GW}^Q(T, \eta)$ and $\mathcal{GW}^Q(T^{\prime}, \eta^{\prime})$ are homotopic.
\end{Prop}

\begin{proof}
The proofs of the corresponding results in the exact setting \cite[\S 2.8]{schlichting2010b} carry over.
\end{proof}

An obvious modification of Lemma \ref{lem:reflectExactIso} (and its proof) is as follows.

\begin{Lem}
\label{lem:reflectExactIsoDual}
Let $\mathcal{A}$ and $\mathcal{B}$ be proto-exact categories with duality and $(F,\eta): \mathcal{A} \rightarrow \mathcal{B}$ an essentially surjective proto-exact form functor which is bijective on inflations and deflations. Then the induced map $\mathcal{GW}^Q(F, \eta) : \mathcal{GW}^Q(\mathcal{A}) \rightarrow \mathcal{GW}^Q(\mathcal{B})$ is a homotopy equivalence.
\end{Lem}

Suppose now that $(\mathcal{A},\oplus)$ is a symmetric monoidal category with duality. Orthogonal direct sum gives $\mathcal{A}_h$ the structure of a symmetric monoidal groupoid. As in \cite[\S 2]{hornbostel2002}, the direct sum Grothendieck--Witt theory space is the group completion
\[
\mathcal{GW}^{\oplus}(\mathcal{A})
=B(\mathcal{A}_h^{-1} \mathcal{A}_h),
\]
with associated Grothendieck--Witt and Witt groups $GW_i^{\oplus}(\mathcal{A}) = \pi_i \mathcal{GW}^{\oplus}(\mathcal{A})$ and $W_i^{\oplus}(\mathcal{A}) = \coker(K^{\oplus}_i(\mathcal{A}) \xrightarrow[]{H} GW_i^{\oplus}(\mathcal{A}))$, $i \geq 0$, respectively. Note that these are indeed groups. We remark that the obvious analogues of Propositions \ref{prop:GWHyper}
and \ref{prop:GWFunctoriality} and Lemma \ref{lem:reflectExactIsoDual} hold for direct sum Grothendieck--Witt theory.

Let $Q_H(\mathcal{A}) \subset Q_h(\mathcal{A})$
be the full subcategory on hyperbolic objects and $\mathcal{GW}^Q_H(\mathcal{A})$ the homotopy fibre of $B Q_H(\mathcal{A}) \rightarrow BQ(\mathcal{A})$ over $0$. Let also
\[
\mathcal{GW}_H^{\oplus}(\mathcal{A}) = B(\mathcal{A}_H^{-1} \mathcal{A}_H),
\]
where $\mathcal{S}_H$ is the symmetric monoidal groupoid of hyperbolic symmetric forms. The following result plays the role of the Group Completion Theorem for the Grothendieck--Witt theory of uniquely split proto-exact categories. Compare with \cite[Theorem 4.2]{schlichting2004}, \cite[Theorem A.1]{schlichting2017} and \cite[Theorem 6.6]{schlichting2019} in the split exact setting.

\begin{Thm}[{\cite[Theorems 3.2 and 3.11]{eberhardtLorscheidYoung2020a}}]
\label{thm:GWComputation}
Let $\mathcal{A}$ be a uniquely split proto-exact category with duality
\begin{enumerate}[wide,labelwidth=!, labelindent=0pt,label=(\roman*)]
\item There is a weak homotopy equivalence $\mathcal{GW}^Q_H(\mathcal{A}) \simeq \mathcal{GW}_H^{\oplus}(\mathcal{A})$.

\item If, moreover, $\mathcal{A}$ is combinatorial and noetherian, then there is a weak homotopy equivalence
\[
\mathcal{GW}^Q(\mathcal{A})
\simeq
\bigsqcup_{w \in W^Q_0(\mathcal{A})}BG_{S_w}\times\mathcal{GW}^Q_H(\mathcal{A}),
\]
where $S_w$ is an isotropically simple representative of the Witt class $w \in W^Q_0(\mathcal{A})$ with self-isometry group $G_{S_w}$.
\end{enumerate}
\end{Thm}


\section{$K$-theory and Grothendieck--Witt theory of pointed monoids}
\label{sec:KGWMonoids}

In this section we study the $K$-theory and Grothendieck--Witt theory of proto-exact categories of projective modules over pointed monoids. The specialization of this section to commutative pointed monoids is the local model for the scheme theoretic considerations of Sections \ref{sec:KThySch} and \ref{sec:GWThySch}.

\subsection{Pointed monoids and their module categories}
\label{sec:monoids}

We record basic material about pointed monoids and their module categories. A detailed reference for commutative pointed monoids is \cite[\S 2]{chu2012}. Many of the results of \cite{chu2012}, and their proofs, apply with only minor changes in the non-commutative setting. See also \cite{knauer1971}, which treats modules over non-commutative semigroups with identity.

A pointed monoid is a semigroup $A$ with a zero (or absorbing element) $0$ and an identity $1$, so that $a \cdot 0 = 0 = 0 \cdot a$ and $a \cdot 1 = a = 1 \cdot a$ for all $a \in A$. A homomorphism of pointed monoids is a semigroup homomorphism which preserves the zero and identity. Pointed monoids and their homomorphisms form a category $\widetilde{\mathcal{M}}_0$. The full subcategory of commutative pointed monoids is $\mathcal{M}_0$.

Let $I$ be an ideal of a commutative pointed monoid $A$, that is, $0 \in I$ and $I A = I$. The quotient pointed monoid $A \slash I$ is the set $(A \setminus I) \cup \{0\}$ with the multiplication $a \cdot b = ab$ if $a,b, ab \in A \setminus I$ and $a \cdot b = 0$ otherwise.

A left $A$-module (also called an $A$-set) is a pointed set $M$, with basepoint denoted again by $0$, together with a left $A$-action under which $0 \in A$ and $1 \in A$ act by the zero and identity map of $M$, respectively. Right $A$-modules are defined similarly. Unless mentioned otherwise, by an $A$-module we mean a left $A$-module. An $A$-module homomorphism is a pointed $A$-equivariant map. Let $A \Modu$ be the category of left $A$-modules and their homomorphisms and $A\modu$ its full subcategory of finitely generated $A$-modules. An $A$-module $P$ is called projective if, for every $A$-module homomorphism $f: P \rightarrow M$ and  surjective $A$-module homomorphism $g: N \rightarrow M$, there exists an $A$-module homomorphism $h: P \rightarrow N$ satisfying $g \circ h = f$. Let $A \proj \subset A \Modu$ be the full subcategory of finitely generated projective $A$-modules.

\begin{Lem}[{\cite[Proposition 2.27]{chu2012}}]
\label{lem:ProjMod}
Every projective $A$-module is of the form $\bigoplus_{i\in J} Ae_i$ where $e_i^2=e_i$ are idempotents in $A$.
\end{Lem}

An $A$-module homomorphism $f: M \rightarrow N$ is called normal if $f^{-1}(n)$ is empty or a singleton for each $n \in N \setminus \{0\}$. This definition of normality is compatible with the categorical definition in the case of monomorphisms and epimorphisms; cf.\ \cite[Proposition 2.15]{chu2012}.
The zero and identity morphisms are normal, as are compositions of normal morphisms. Denote by $A \Modu^n \subset A \Modu$ the subcategory of normal $A$-module homomorphisms, and similarly for $A\modu^{\nor}$ and $A\proj^{\nor}$.

An element $a \in A$ is called right cancellative (resp. right partially cancellative, or rpc) if right multiplication $\cdot a: A\rightarrow A$ is an injective $A$-module homomorphism (resp. normal $A$-module homomorphism). Explicitly, $a \in A$ is $\rpc$ if $xa=ya$ implies $x=y$ or $xa=ya=0$ for all $x,y\in A$. Let $A^{\rpc} \subset A$ be the subset of right partially cancellative elements and $A^{\rc} \subset A$ be the subset of right cancellative elements together with $0 \in A$. Both $A^{\rc}$ and $A^{\rpc}$ are pointed submonoids of $A$. We call a pointed monoid $A$ right cancellative if $A^{\rc} = A$ and rpc if $A^{\rpc}=A$. Replacing right with left multiplication leads to the notion of a left (partially) cancellative pointed monoid. A pointed monoid is called cancellative (resp. partially cancellative, or simply pc) if it is both left and right cancellative (resp. partially cancellative).

A pointed monoid $A$ is called right reversible
if $Aa\cap Ab\neq \{0\}$ for any two non-zero elements $a,b\in A$. For example, a commutative cancellative pointed monoid is (both left and right) reversible, since $0 \neq ab \in Aa \cap A b$.

A pointed monoid is called right noetherian if it satisfies the ascending chain condition for right congruences.

For a family of $A$-modules $\{M_i\}_{i\in J}$, the direct sum $A$-module is
\[
\bigoplus_{i\in J} M_i = \Big( \bigsqcup_{i\in J} M_i  \Big) \slash \langle 0_{M_i} \sim 0_{M_j} \mid i,j \in J \rangle
\]
with the obvious $A$-action. For an $A\mhyphen B$-bimodule $M$ and a $B \mhyphen C$-bimodule $N$, the tensor product $A \mhyphen C$-bimodule is
\[
M \otimes_B N = \left( M \times N \right) \slash \{ (mb,n) \sim (m, bn) \mid b \in B \}
\]
with the obvious actions of $A$ and $C$.

\begin{Ex}
\begin{enumerate}[wide,labelwidth=!, labelindent=0pt,label=(\roman*)]
\item The initial object of $\widetilde{\mathcal{M}}_0$ is $\mathbb{F}_1:=\{0,1\}$. There is an equivalence of $\Fun \modu$ with the category $\mathsf{set}_*$ of finite pointed sets. The subcategory $\Fun \modu^{\nor} = \Fun \proj^{\nor}$ is often denoted by $\Vect_{\Fun}$ in the literature.

\item The terminal object of $\widetilde{\mathcal{M}}_0$ is $\{0\}$, the unique monoid with $0=1$.

\item Let $\mathsf{G}$ be a group. Then $\Fun[\mathsf{G}] := \mathsf{G} \sqcup \{0 \}$ is a cancellative pointed monoid. A pointed monoid is cancellative and right reversible if and only if it can be embedded in $\Fun[\mathsf{G}]$ for some group $\mathsf{G}$ \cite[Theorem 1.23]{clifford1961}.

\item The subset $A^{\times} \subset A$ of multiplicative units is a group and $\Fun[A^{\times}] \subset A$ is a cancellative pointed submonoid.

\item The pointed monoid $\Fun[t] = \{t^i\}_{i \geq 0} \sqcup \{0\}$ is cancellative.

\item Let $n \geq 2$. The pointed monoid $\Fun[t] \slash \langle t^n = 0 \rangle = \{ 0,1,t, \dots, t^{n-1} \}$
is not cancellative, since $t \cdot 0 = t \cdot t^{n-1}$, but is pc and reversible.

\item The pointed monoid $A= \Fun[t,s] \slash \langle ts =0\rangle$ is pc but not reversible, since $A t \cap A s = \{0\}$.

\item Let $n > d \geq 2$. The pointed monoid $\Fun[t] \slash \langle t^n = t^d \rangle = \{ 0,1,t,\dots, t^d, \dots, t^{n-1} \}$ is not pc, since $t \cdot t^{d-1} = t \cdot t^{n-1}$.
\end{enumerate}
\end{Ex}

For a left $A$-module $M$, the set $\Hom_{A \Modu}(M,A)$ becomes a right $A$-module via
\[
(f\cdot a)(m):=f(m)a,
\qquad
f\in \Hom_{A \Modu}(M,A), \; a\in A, \; m\in M.
\]
Unlike in the case of rings, the module $\Hom_{A \Modu}(M,A)$ does not define a good notion of a module dual to $M$. For this reason, we instead consider the subset $\Hom_{A \Modu^\nor}(M,A) \subset \Hom_{A \Modu}(M,A)$ of normal homomorphisms. As the following result shows, this subset is not an $A$-submodule without additional assumptions. Denote by $A^{\op}$ the monoid opposite to $A$.

\begin{Lem}
Let $A$ be a pointed monoid.
\label{lem:dualFreeAMod}
\begin{enumerate}[wide,labelwidth=!, labelindent=0pt,label=(\roman*)]
\item For any $M \in A \Modu$, the right $A$-module structure on $\Hom_{A \Modu}(M,A)$ induces a right $A^{\rpc}$-module structure on $\Hom_{A \Modu^\nor}(M,A)$.

\item If $A$ is right reversible and rpc and $M$ is a finitely generated free $A$-module, then the right $A$-module $\Hom_{A \modu^{\nor}}(M,A)$ is finitely generated and free.

\item \label{part:oplusMonoidal} If $A$ is right reversible, rpc and right noetherian, 
then $\Hom_{A \modu^{\nor}}(-, A)$ defines a $\oplus$-monoidal functor
\[
P : (A \modu^{\nor})^{\op} \rightarrow A^{\op} \modu^{\nor}.
\]
\end{enumerate}
\end{Lem}

\begin{proof}
The first statement is a direct verification. For the second statement, let
$$M=\bigoplus_{i\in J} As_i$$
be a finitely generated free $A$-module. Denote by $s_i^{\vee}:M\rightarrow A$ the map sending $as_i$ to $a$ and $as_j$ to $0$ if $j\neq i$. We claim that the induced map
\[
\bigoplus_{i\in J}s_i^{\vee}: \bigoplus_{i\in J}t_iA \rightarrow \Hom_{A \modu^{\nor}}(M,A),
\qquad
t_i x_i \mapsto s_i^{\vee} \cdot x_i
\]
is a right $A$-module isomorphism. The map is well-defined and injective since $A$ is rpc. If $|J|=1$, then the map is clearly an isomorphism. Suppose then that $|J|\geq 2$ and let $f\in \Hom_{A \Modu^{\nor}}(M,A)$. Set $f_i=f(s_i)$. We claim that there is at most one $i\in J$ such that $f_i\neq 0$ and hence $f=s_i^{\vee}\cdot f_i$. Assume that there exist distinct $i, j \in J$ such that $f_i\neq 0\neq f_j$. Since $A$ is right reversible, there exist $a,b\in A$ such that $af_i=bf_j\neq 0$. Hence, $f(as_i)=f(bs_j)\neq 0$, a contradiction. The second statement follows.

Turning to the third statement, let $M \in A \modu$. Fix a surjection $F \rightarrow M$ with $F$ a finitely generated free $A$-module. A direct check shows that $P(M)$ is naturally a submodule of $P(F)$. By the first two parts of the lemma, $P(F)$ is finitely generated and free. Since $A$ is right noetherian, $P(F)$ is noetherian \cite[Proposition 2.31]{chu2012}, from which it follows that $P(M)$ is finitely generated. The definition of $P$ on morphisms is via pre-composition and is well-defined because the composition of normal morphisms is normal. To prove that $P$ is $\oplus$-monoidal, let $M, N \in A \modu$. An element $f \in P(M \oplus N)$ determines by restriction $f_M \in P(M)$ and $f_N \in P(N)$. Suppose that neither $f_M$ nor $f_N$ is zero. Since $A$ is right reversible, $\im \, f_M \cap \im \, f_N \neq \{0\}$, contradicting the assumption that $f$ is normal. It follows that at most one of $f_N$ and $f_M$ is non-zero and there is a well-defined $A$-module homomorphism
\[
P(M \oplus N) \rightarrow
P(M) \oplus
P(N).
\]
It is straightforward to verify that this is an isomorphism. We omit the verification that $P$ respects $\oplus$ on morphisms.
\end{proof}

\begin{Rem}
\begin{enumerate}[wide,labelwidth=!, labelindent=0pt,label=(\roman*)]
\item There is a right $A$-module isomorphism
\[
\Hom_{A \modu}(\bigoplus_{i \in J} A s_i,A) \simeq \prod_{i \in J} A.
\]
In particular, the standard $A$-linear dual of a free $A$-module is in general not free. In fact, $\prod_{i \in J} A$ need not even be finitely generated. For example, the $\Fun[t]$-module $\Fun[t] \times \Fun[t]$ is not finitely generated. From this point of view, the normal dual $\Hom_{A \modu^{\nor}}(-, A)$ has better properties than $\Hom_{A \modu}(-,A)$.

\item The functor $\Hom_{A \modu^{\nor}}(-, A)$ of Lemma \ref{lem:dualFreeAMod}\ref{part:oplusMonoidal} does not extend to $(A \modu)^{\op} \rightarrow A^{\op} \modu$, since for a non-normal morphism $f: M \rightarrow N$, the image of $\Hom_{A \modu^{\nor}}(f, A)$ is not contained in $\Hom_{A \modu^{\nor}}(M, A) \subset \Hom_{A \modu}(M,A)$. 
\end{enumerate}
\end{Rem}

\subsection{$K$-theory of pointed monoids}
\label{sec:KThyMonoids}

The $K$-theory of pointed monoids has been studied by a number of authors \cite{deitmar2006}, \cite{chu2012}, \cite{chu2010}, \cite{haesemeyer2019}. In this section we describe those results which are relevant to this paper.

Let $A$ be a pointed monoid. The category $A \Modu$ admits a proto-exact structure with inflations and deflations being the normal $A$-module homomorphisms which are injective and surjective, respectively \cite[\S 2.2.2]{chu2012}.
 However, $\oplus$ is not a coproduct for $A \Modu^{\nor}$. Indeed, for a non-zero $A$-module $M$, there is no dashed arrow in $A \Modu^{\nor}$ which makes the diagram
\[
\begin{tikzpicture}[baseline= (a).base]
\node[scale=1.0] (a) at (0,0){
\begin{tikzcd}[column sep=4ex,row sep=4ex]
& M \arrow[tail,"i_M" left]{d} \arrow[bend right=-20]{rdd}{\id_M}& \\
M \arrow[tail, "i_M"]{r} \arrow[bend right=20,swap]{rrd}{\id_M} & M \oplus M \arrow[dashed]{rd}& \\
& & M
\end{tikzcd}
};
\end{tikzpicture}
\]
commute. In particular, $A \Modu^{\nor}$ is not a quasi-exact category.

Since $A \proj \subset A \Modu$ is an extension closed full subcategory, it inherits a proto-exact structure from $A \Modu$.

\begin{Lem}[{See also \cite[Proposition 2.29]{chu2012}.}]
\label{lem:AProjSplit}
Let $A$ be a pointed monoid. The proto-exact category $A \proj$ is uniquely split and combinatorial.
\end{Lem}

\begin{proof}
Let $U \xrightarrowtail[]{i} V \xtwoheadrightarrow[]{\pi} W$ be a conflation in $A \proj$. Since $W$ is projective, there exists a section $s: W \rightarrow V$ of $\pi$. Since $\pi$ is a deflation, it is normal by definition which implies that that the section $s$ is unique. Define an $A$-module homomorphism $\phi: U \oplus W \rightarrow V$ by
\[
\phi(u) = i(u), \qquad \phi(w) = s(w).
\]
To see that $\phi$ is injective, suppose, for example, that $\phi(u) = \phi(w)$. Applying $\pi$ gives
\[
0=\pi(i(u))=\pi(s(w)) = w,
\]
implying $u=w=0$. We claim that $\phi$ is also surjective and hence an isomorphism by \cite[Lemma 2.2]{chu2012}. It is immediate that $\im \, i \subset \im \, \phi$. Let $v \in V \setminus \im \, i$. Then $\pi(v) \neq 0$ and hence also $s(\pi(v)) \neq 0$. Because $\pi \circ s = \id_W$, the map $s \circ \pi$ is idempotent. It follows that $s(\pi(v))$ and $v$ have the same (non-zero) image under $s\circ \pi$. Since $s \circ \pi$ is normal, $s(\pi(v))=v$. We conclude that $\phi$ is a splitting of the original conflation.

That the combinatorial property holds follows from the fact that $\oplus$ is defined using disjoint union of the underlying sets.
\end{proof}

Note that Lemma \ref{lem:AProjSplit} also implies that $A \proj^{\nor}$ is a uniquely split proto-exact category. The following `$Q=+$' theorem is the main results of this section.

\begin{Thm}
\label{thm:KThyRPC}
Let $A$ be an rpc pointed monoid. Then there is a homotopy equivalence
\[
\mathcal{K}(A \proj)
\simeq
\mathbb{Z} \times B (A^{\times} \wr \Sigma_{\infty})^+.
\]
In particular, if $A^{\times}$ is finite, then $K_0(A \proj) \simeq \mathbb{Z}$ and
\[
K_i(A \proj) \simeq \pi_i^s(B A^{\times}_+),
\qquad
i \geq 1.
\]
\end{Thm}

\begin{proof}
Since $A$ is $\rpc$, projective $A$-modules are free. Indeed, this follows from Lemma \ref{lem:ProjMod} and the fact that a (non-trivial) $\rpc$ pointed monoid has a single non-zero idempotent, namely $1 \in A$. By Theorem \ref{thm:KComparison} (see also \cite[Theorem 4.2]{chu2010}), there is a homotopy equivalence $\mathcal{K}(A \proj) \simeq \mathcal{K}^{\oplus}(A \proj)$. To compute $\mathcal{K}^{\oplus}(A \proj)$, we apply \cite[Proposition 3]{weibel1981} to the cofinal family $\{A^{\oplus n}\}_{n \in \mathbb{Z}_{\geq 0}}$ of $A \proj$. We then have
\[
\Aut(A \proj) : = \varinjlim_n \, \Aut_{A \proj}(A^{\oplus n}) = \varinjlim_n \, (A^{\times})^n \rtimes \Sigma_n =  A^{\times} \wr \Sigma_{\infty},
\]
giving the claimed result.
\end{proof}

Note that, by Lemma \ref{lem:reflectExactIso}, the embedding $A \proj^{\nor} \hookrightarrow A \proj$ induces a homotopy equivalence $\mathcal{K}(A \proj^{\nor}) \simeq \mathcal{K}(A \proj)$.

\subsection{Grothendieck--Witt theory of pointed monoids}
\label{sec:GWThyMonoids}

In this section, we study the Grothendieck--Witt theory of pointed monoids. This leads to a non-additive analogue of Karoubi's Grothendieck--Witt theory of rings \cite{karoubi1973}.

Let $A$ be an $\rpc$ pointed monoid. Fix a pointed monoid involution $\sigma: A \rightarrow A^{\op}$ and a central element $\epsilon \in A$ which satisfies $\epsilon \sigma(\epsilon) =1$. For example, when $A$ is commutative, $\sigma = \id_A$ and $\epsilon =1$ is an admissible choice. For a non-trivial case, see the examples below.

Given $M \in A \Modu$, consider $P^{\sigma}(M) := \Hom_{A \Modu^n}(M,A)$ as a left $A$-module via
\[
(a \cdot f)(m) := f(m) \sigma(a),
\qquad
f \in P^{\sigma}(M), \, a \in A, \, m \in M.
\]
Compare with Lemma \ref{lem:dualFreeAMod}.

\begin{Prop}
\label{prop:projDuality}
Let $A$ be a right reversible rpc pointed monoid. The natural transformation $\Theta^{\sigma, \epsilon}: \id_{A \proj^{\nor}} \Rightarrow P^{\sigma} \circ (P^{\sigma})^{\op}$ with components
\[
\Theta^{\sigma, \epsilon}_M(m)(f) = \epsilon \sigma(f(m)),
\qquad
f \in P^{\sigma}(M), \, m \in M
\]
makes $(A \proj^{\nor}, P^{\sigma},\Theta^{\sigma, \epsilon})$ into a uniquely split combinatorial proto-exact category with duality.
\end{Prop}

\begin{proof}
Lemma \ref{lem:AProjSplit} shows that $A \proj$, and hence $A \proj^{\nor}$, is uniquely split and combinatorial. Given $\phi: M \rightarrow N$ in $A \proj^{\nor}$, the morphism $P^{\sigma}(\phi): P^{\sigma}(N) \rightarrow P^{\sigma}(M)$ is defined to be $(-) \circ \phi$. This is well-defined since the composition of normal morphisms is normal. Since $A$ is $\rpc$, projective $A$-modules are free. That $P^{\sigma}$ is $\oplus$-monoidal on $A \proj^{\nor}$ follows from Lemma \ref{lem:dualFreeAMod}. Exactness of $P^{\sigma}$ then follows from the splitness of $A \proj^{\nor}$. Hence, $P^{\sigma}$ satisfies the desired properties.

A direct calculation shows that $\Theta^{\sigma,\epsilon}_M(m) : P^{\sigma}(M) \rightarrow A$ is an $A$-module homomorphism. To see that $\Theta^{\sigma, \epsilon}_M(m)$ is normal, fix an $A$-module basis $M \simeq \bigoplus_{i \in J} A s_i$ and write $m = x s_i$. When $x=0$, the map $\Theta^{\sigma,\epsilon}_M(m)$ is zero, which is normal. Suppose then that $x \neq 0$ and let $0 \neq a \in A$. We have
\begin{multline*}
\Theta^{\sigma, \epsilon}_M(m)^{-1}(a)
=
\{ f \in P^{\sigma}(M) \mid \epsilon(\sigma(f(x s_i)) = a\}
=
\{ y \in A \mid \epsilon(\sigma(y s_i^{\vee}(x s_i)) = a\}
\\
=\{ y \in A \mid \sigma(y) \sigma(x) = \sigma(\epsilon) a\}.
\end{multline*}
Since $A$ is $\rpc$ and $\sigma$ is an isomorphism, the final set is empty or a singleton, as required. The assumption $\epsilon \sigma(\epsilon)=1$ ensures that the equalities $P(\Theta_U) \circ \Theta_{P(U)} = \id_{P(U)}$ hold.
\end{proof}

\begin{Ex}
Let $A = \Fun$. The only possibilities are $\sigma=\id_{\Fun}$ and $\epsilon =1$. For each $M \in \Vect_{\Fun}=\Fun \proj^{\nor}$, there is a canonical isomorphism
\[
\delta_M: M \xrightarrow[]{\sim} P(M),
\qquad
\delta_M(m)(m^{\prime})  = 
\begin{cases}
1 & \mbox{ if } m= m^{\prime}, \\
0 & \mbox{ if } m\neq m^{\prime}.
\end{cases}
\]
We emphasize that such an isomorphism does not exist for a general pointed monoid. Under this identification, $P$ squares to the identity. The triple $(\Vect_{\Fun}, P, \id_{\id_{\Vect_{\Fun}}})$ is therefore a proto-exact category with strict duality.
\end{Ex}

Next, we turn to the classification of symmetric forms in $A \proj^{\nor}$. Let $M$ be a free $A$-module of rank one. Fix a basis $M \simeq A$. A symmetric form $\psi_M$ on $M$ then takes the form
\[
\psi_M(a)(x) = x \xi \sigma(a),
\qquad
a, x \in A
\]
for some $\xi \in A^{\times}$ which satisfies $\xi = \epsilon \sigma(\xi)$. Write $\psi_{\xi}$ for this symmetric form. An isomorphism $M \rightarrow M$, which is necessarily determined by an element $u \in A^{\times}$, defines an isometry $\psi_{u \xi \sigma(u)} \rightarrow \psi_{\xi}$. Motivated by these observations, define an $A^{\times}$-action on the set
\[
A^{\times}_{\sigma, \epsilon} = \{ \xi \in A^{\times} \mid \xi = \epsilon \sigma(\xi)\}
\]
by $u \cdot \xi = u \xi \sigma(u)$. Then the set of isomorphism classes of rank one symmetric forms in $A \proj^{\nor}$ is $\SPic(A): =A^{\times}_{\sigma, \epsilon} \slash A^{\times}$. The isometry group of $\psi_{\xi}$ is the stabilizer
\[
I(\xi) = \{ u \in A^{\times} \mid \xi = u \xi \sigma(u)\}.
\]

\begin{Ex}
The pointed monoid $A = \mathbb{F}_{1^3} := \Fun[\mathbb{Z} \slash 3]$ has a unique non-trivial monoid automorphism $\sigma$, which is an involution. Either non-identity element $\epsilon \in \mathbb{Z} \slash 3$ is compatible with $\sigma$. We have $A^{\times}_{\sigma, \epsilon} = \{\epsilon^2\}$ and $I (\epsilon^2) \simeq \mathbb{Z} \slash 3$. In particular, $\SPic(A)$ is a singleton.
\end{Ex}

Given $(h,\{m_{\xi}\}) \in \mathbb{Z}_{\geq 0} \times \sideset{}{'}\prod_{\xi \in \SPic(A)} \mathbb{Z}_{\geq 0}$ in the restricted product (see Appendix \ref{app:RestrictedProduct}), define a symmetric form
\[
\psi_{h,\{m_{\xi}\}}
=
\psi_{H(A)}^{\oplus h} \oplus
\bigoplus_{\xi \in \SPic(A)} \psi_{\xi}^{\oplus m_{\xi}}.
\]

\begin{Prop}
\label{prop:autGroupAProj}
Let $A$ be a right reversible rpc pointed monoid.
\begin{enumerate}[wide,labelwidth=!, labelindent=0pt,label=(\roman*)]
\item The assignment $(h,\{m_{\xi}\}) \mapsto \psi_{h,\{m_{\xi}\}}$ induces a monoid isomorphism between
\[
\mathbb{Z}_{\geq 0} \times \sideset{}{'}\prod_{\xi \in \SPic(A)} \mathbb{Z}_{\geq 0}
\]
and the monoid $\pi_0(A \proj^{\nor}_h)$ of isometry classes of symmetric forms in $A \proj^{\nor}$.

\item There is a group isomorphism
\[
\Aut_{A \proj^{\nor}_h}(\psi_{h,\{m_\xi\}})
\simeq
\big(
(\mathbb{Z} \slash 2 \ltimes_{\sigma} A^{\times}) \wr \Sigma_h
\big)
\times
\sideset{}{'}\prod_{\xi \in \SPic(A)} \big( I(\xi)  \wr \Sigma_{m_{\xi}} \big),
\]
where $\mathbb{Z} \slash 2$ acts on $A^{\times}$ by $u \mapsto \sigma(u^{-1})$.
\end{enumerate}
\end{Prop}

\begin{proof}
After using Lemma \ref{lem:noMeta}, the first statement is straightforward. The second statement is a direct calculation.
\end{proof}

\begin{Rem}
Let $M \in A \proj$ be free of rank $n$. Fixing a basis of $M$, and hence also of $P^{\sigma}(M)$, identifies a symmetric form on $M$ with an $A^{\times}$-valued permutation matrix $\psi=(\psi_{ij}) \in A^{\times} \wr \Sigma_n$ which satisfies $\psi_{ij} = \epsilon \sigma(\psi_{ji})$, $1 \leq i,j \leq n$. In this formulation, Proposition \ref{prop:autGroupAProj} becomes the classification of such matrices up to congruence.
\end{Rem}

Recall the Reduction Assumption from Section \ref{sec:catWD}.

\begin{Prop}
\label{prop:redAssumptionMonoid}
The Reduction Assumption holds for $(A \proj^{\nor}, P^{\sigma},\Theta^{\sigma, \epsilon})$.
\end{Prop}

\begin{proof}
An isotropic subobject $U \rightarrowtail \psi_{h,\{m_{\xi}\}}$ necessarily factors through the summand $H(A)^{\oplus h}$. It therefore suffices to consider only hyperbolic symmetric forms. By the combinatorial property of $A \proj^{\nor}$, we can write $U = U_1 \oplus P(U_2)$ for some $U_i \in A \proj$, in which case the isotropic condition is $P(U_2) \rightarrowtail P(X \slash U_1)$. The reduction of $H(A)^{\oplus h}$ is then canonically isometric to $H(\coker (P(U_2) \rightarrowtail P(X \slash U_1))$. See also \cite[Lemma 1.1]{mbyoung2021}.
\end{proof}

\begin{Rem}
In view of Proposition \ref{prop:redAssumptionMonoid}, we conclude, using \cite[Theorem 3.10]{mbyoung2018b}, that the forgetful morphism $\mathcal{R}_{\bullet}(A \proj^{\nor}) \rightarrow \mathcal{S}_{\bullet}(A \proj^{\nor})$ from the $\mathcal{R}_{\bullet}$-construction to the Waldhausen $\mathcal{S}$-construction is a relative $2$-Segal space. We can therefore apply the construction of \cite[\S 4]{mbyoung2018b} to produce a module over of the Hall algebra $\mathcal{H}(A \proj^{\nor})$. The algebra $\mathcal{H}(A \proj^{\nor})$, and its variations, have been studied by Szczesny \cite{szczesny2012}, \cite{szczesny2014}, \cite{szczesny2018}. In the setting of the representation theory of quivers over $\mathbb{F}_1$, which is combinatorial but not split, modules arising from the $\mathcal{R}_{\bullet}$-construction have been studied in \cite{mbyoung2021} where, in particular, a version of Green's theorem is proved. \end{Rem}

We can now state the main result of this section.

\begin{Thm}
\label{thm:GWThyMonoid}
Let $A$ be a right reversible rpc pointed monoid with $\SPic(A)$ countable. Then $\mathcal{GW}^{\oplus}(A \proj^{\nor}, \sigma, \epsilon)$ is homotopy equivalent to
\[
\mathbb{Z} \times \Big(\sideset{}{'}\prod_{\xi \in \SPic(A)} \mathbb{Z} \Big) \times B\left( \big( (\mathbb{Z} \slash 2 \ltimes_{\sigma} A^{\times}) \wr \Sigma_{\infty} \big) \times \sideset{}{'}\prod_{\xi \in \SPic(A)}  \big(I(\xi) \wr \Sigma_{\infty} \big) \right)^+.
\]
\end{Thm}

\begin{proof}
Index the set $\SPic(A)$ as $\{\xi_j\}_{j\in J}$ for some subset $J\subseteq \mathbb{Z}_{\geq 0}$. For each $n\geq 0$, consider the symmetric form
\[
\Psi_n = \bigoplus_{\substack{j\in J\\ j\leq n}} \psi_{\xi_j}
\]
and set $s_n = H(A)^{\oplus n} \oplus \Psi_n^{\oplus n}$. Then $\{s_n\}_{n \in \mathbb{Z}_{\geq 0}}$ is a cofinal family in $A \proj^{\nor}_h$. 
Using Proposition \ref{prop:autGroupAProj}, we find
\[
\Aut(A \proj^{\nor}_h) := \varinjlim \Aut_{A \proj^{\nor}_h}(s_n)
\simeq
\big(
(\mathbb{Z} \slash 2 \ltimes_{\sigma} A^{\times}) \wr \Sigma_{\infty}
\big)
\times \sideset{}{'}\prod_{\xi \in \SPic(A)}  I(\xi) \wr \Sigma_{\infty} .
\]
We are therefore in the setting of \cite[Proposition 3]{weibel1981}, allowing us to conclude that there is a homotopy equivalence
\[
\mathcal{GW}^{\oplus}(A \proj^{\nor}, \sigma, \epsilon)
\simeq K_0^{\oplus}(A \proj^{\nor}_h) \times B\Aut(A \proj^{\nor}_h)^+.
\]
Finally, use Proposition \ref{prop:autGroupAProj} to identify $K_0^{\oplus}(A \proj^{\nor}_h)$ and $\mathbb{Z} \times \sideset{}{'}\prod_{\xi \in \SPic(A)} \mathbb{Z}$.
\end{proof}

\begin{Cor}
\label{thm:WittThyMonoid}
In the setting of Theorem \ref{thm:GWThyMonoid}, there is an isomorphism
\[
W^{\oplus}_0(A \proj^{\nor}, \sigma, \epsilon)
\simeq
\sideset{}{'}\prod_{\xi \in \SPic(A)} \mathbb{Z}.
\]
\end{Cor}

\begin{Ex}
Suppose that $A$ has no non-trivial units. This is the case, for example, for $A = \Fun[T_1, \dots, T_n]$. Then there is an isomorphism
\[
\mathcal{GW}^{\oplus}(A \proj^{\nor}, \sigma)
\simeq
\mathbb{Z}^2 \times B\left( \big( \mathbb{Z} \slash 2 \wr \Sigma_{\infty} \big) \times \Sigma_{\infty} \right)^+.
\]
\end{Ex}

\begin{Ex}
For later use, we record the homotopy equivalence
\[
\mathcal{GW}^{\oplus}(\Vect_{\mathbb{F}_1})
\simeq
\mathbb{Z}^2 \times B\left( \big( \mathbb{Z} \slash 2 \wr \Sigma_{\infty} \big) \times \Sigma_{\infty} \right)^+.
\]

To describe $\mathcal{GW}^Q(\Vect_{\mathbb{F}_1})$, we use Theorem \ref{thm:GWComputation}. Since $\{H(\mathbb{F}_1^{\oplus n})\}_{n \geq 0}$ is a cofinal family of $\mathcal{A}_H$, arguing as in the proof of Theorem \ref{thm:GWThyMonoid}, we obtain a homotopy equivalence
\[
\mathcal{GW}_H(\Vect_{\mathbb{F}_1})
\simeq
\mathbb{Z} \times B (\mathbb{Z} \slash 2 \wr \Sigma_{\infty})^+
\]
so that there is a weak homotopy equivalence
\[
\mathcal{GW}^Q(\Vect_{\mathbb{F}_1})
\simeq
\bigsqcup_{n \in \mathbb{Z}_{\geq 0}}B \Sigma_n \times \mathbb{Z} \times B (\mathbb{Z} \slash 2 \wr \Sigma_{\infty})^+.
\]
Computations for more general monoids are similar.
\end{Ex}


\section{Algebraic $K$-theory of monoid schemes}
\label{sec:KThySch}

In the remainder of the paper, all pointed monoids are assumed to be commutative.
 
\subsection{Monoid schemes}
\label{sec:F1Sch}

We present some background on monoid schemes. The reader is referred to \cite{chu2012}, \cite{connes2010}, \cite{cortinas2015}, \cite{deitmar2005} for further details.

A prime ideal of a pointed monoid $A$ is an ideal $\fp \subset A$ whose complement $S=A-\fp$ is a multiplicative subset, that is, $S$ contains $1$ and is multiplicatively closed.

Let $S$ be a multiplicative subset of $A$. The localization of $A$ at $S$ is $S^{-1}A=(S\times A)/\sim$, where $(s,a)\sim(a',s')$ if there exists $t\in S$ such that $tsa'=ts'a$. Write $\frac as$ for the class of $(s,a)$ in $S^{-1}A$. The product $\frac as\cdot\frac bt=\frac{ab}{st}$ endows $S^{-1}A$ with the structure of a pointed monoid and the map $\iota_S:A\to S^{-1}A$, $a \mapsto \frac{a}{1}$, is a monoid morphism, which we call the localization map.

Given $h\in A$, we write $A[h^{-1}]$ for the localization of $A$ at $S=\{h^i\}_{i\in \mathbb{Z}_{\geq 0}}$. Given a prime ideal $\fp \subset A$, we write $A_\fp$ for the localization of $A$ at $S=A-\fp$. If $A$ is cancellative, then $S=A-\{0\}$ is a multiplicative subset and we define the fraction field of $A$ as $\Frac A=S^{-1}A$, which is a pointed group, 
that is, $(\Frac A)^\times=\Frac A-\{0\}$. More generally, if $A$ is cancellative and $0\notin S$, then $\iota_S:A\to S^{-1}A$ is injective and $S^{-1}A$ is cancellative. In this situation, we often identify $A$ with its image in $S^{-1}A$ and write $a$ for $\frac a1\in S^{-1}A$.

A monoidal space is a pair $(X,\mathcal{O}_X)$ consisting of a topological space $X$ and a sheaf of pointed monoids $\mathcal{O}_X$. We often surpress $\cO_X$ from the notation. A primary example of a monoidal space is the spectrum $X=\Spec A$ of a pointed monoid $A$ whose points are the prime ideals of $A$, whose topology is generated by the principal open subsets $U_h=\{\fp\mid h\notin\fp\}$ for $h\in A$, and whose structure sheaf $\cO_X$ is characterized by the values $\cO_X(U_h)=A[h^{-1}]$ and by its stalks $\cO_{X,\fp}=A_\fp$. We use the short hand notation $\Gamma X=\cO_X(X)$ for the pointed monoid of global sections of $\mathcal{O}_X$.

An affine monoid scheme is a monoidal space which is isomorphic to $\Spec(A)$ for some pointed monoid $A$. A monoid scheme is a monoidal space which admits an open cover by affine monoid schemes.

A morphism between monoid schemes $X$ and $Y$ is a continuous map $\varphi:X\to Y$ together with a morphism $\varphi^\#:\cO_Y\to\varphi_\ast(\cO_X)$ of sheaves of pointed monoids such that, for every $x\in X$, the induced pointed monoid morphism $\varphi_x^\#:\cO_{Y,\varphi(x)}\to\cO_{X,x}$ is local, that is, maps non-units to non-units.

A monoid scheme is of finite type if it is quasi-compact and has an affine open covering by spectra of finitely generated pointed monoids. We say that $X$ has enough closed points if every point $y \in X$ specializes to a closed point $x$ or, equivalently, $x$ is contained in the topological closure of $\{y\}$.

\begin{Rem}\label{rem:PropsMonSch}
We list some well-known properties of monoid schemes.
\begin{enumerate}[wide,labelwidth=!, labelindent=0pt,label=(\roman*)]
\item \label{lem:irredUniqueGenPt} A monoid scheme is a spectral space (cf.\ \cite{ray2020}), which means, in particular, that every irreducible closed subset has a unique generic point, which we typically denote by $\eta$.

\item A pointed monoid $A$ has a unique maximal (prime) ideal, namely, the complement $\fm=A-A^\times$ of the unit group $A^\times$. Therefore, every affine monoid scheme $X=\Spec A$ has a unique closed point $x=\fm$, and $X$ is the only open neighbourhood of $x$. As a consequence, every affine open subset $U$ of a monoid scheme $X$ has a unique closed point $x$ and $\Gamma U=\cO_{X,x}$. 

\item If $X$ has enough closed points, then it is covered by the minimal open neighbourhoods $U_x=\Spec \cO_{X,x}$ of the closed points $x$. This covering is the minimal open covering of $X$, in the sense that every other open covering of $X$ refines to $\{U_x\mid x \in X \text{ closed}\}$. 

\item If $X$ is of finite type, then it has only finitely many points. In particular, $X$ has enough closed points and $U_x$ is open in $X$ for every $x\in X$.
\end{enumerate}
\end{Rem}

\begin{Def}
A monoid scheme $X$ is called
\begin{enumerate}[wide,labelwidth=!, labelindent=0pt,label=(\roman*)]
\item cancellative (resp. pc) if the pointed monoid $\mathcal{O}_{X,x}$ is cancellative (resp. pc) for each $x \in X$,

\item integral (resp. reversible) if the pointed monoid $\mathcal{O}_X(U)$ is cancellative (resp. reversible) for each open set $U \subset X$, and

\item torsion free if for every $x\in X$, the unit group of $\mathcal{O}_{X,x}$ is torsion free, that is, if $a^n=1$ for $n>1$, then $a=1$.
\end{enumerate}
\end{Def}

\begin{Rem}\label{lem:pcGlobal}\label{rem:CancelMonoidScheme}
Note that we digress from \cite{chu2012} in the meaning of integrality of monoid schemes. To wit, we require that $\cO_X(U)$ is cancellative for all opens $U$, and not merely for an open covering, as required in \cite{chu2012}.
\end{Rem}

\begin{Lem}
\label{lem:pcProperties}
\phantomsection
\begin{enumerate}[wide,labelwidth=!, labelindent=0pt,label=(\roman*)]
\item \label{part:pcLoc} Let $X$ be a cancellative (resp.\ pc) monoid scheme. Then the pointed monoid $\mathcal{O}_X(U)$ is cancellative (resp.\ pc) for each affine open subset $U \subset X$.

\item An integral monoid scheme is cancellative.

\item \label{part:intRever} An integral monoid scheme is reversible.
\end{enumerate}
\end{Lem}

\begin{proof}
Let $U \subset X$ be an open affine subset of a cancellative (resp.\ pc) monoid scheme. Then $U=\Spec\cO_{X,x}$ where $x$ is the unique closed point of $U$. It follows that $\mathcal{O}_X(U)$ is cancellative (resp.\ pc).

The second statement follows from the fact that the stalks of an integral monoid scheme are submonoids of the (cancellative) generic stalk.

The final statement follows from the first statement and the fact that cancellative pointed monoid is reversible.
\end{proof}

Let $X$ be an irreducible cancellative monoid scheme with generic point $\eta$. Define the function field of $X$ as $\cO_{X,\eta}$, which is a pointed group.

\begin{Lem}
\label{lem:genericSubmonoid}
Let $X$ be an irreducible cancellative monoid scheme with generic point $\eta$. For every open $U\subset X$, there is an equality $\cO_X(U)=\bigcap_{x\in U}\cO_{X,x}$ of pointed submonoids of $\cO_{X,\eta}$.
\end{Lem}

\begin{proof}
Since $\cO_{X,x}$ is cancellative and $\cO_{X,\eta}$ is a localization of $\cO_{X,x}$, the localization map $\cO_{X,x}\to \cO_{X,\eta}$ is injective for every $x\in X$.
\end{proof}

\begin{Prop}
\label{prop:intVsCanc}
 A monoid scheme $X$ is integral if and only if it is irreducible and pc.
\end{Prop}

\begin{proof}
Assume that $X$ is integral. If $X$ is not irreducible, then there exist disjoint non-empty open subsets $U_1$ and $U_2$. By the sheaf axiom, we have $\cO_X(U_1\cup U_2)=\Gamma U_1\times\Gamma U_2$, which is not cancellative, contradicting the integrality of $X$. Thus, $X$ is irreducible. Let $x \in X$ with affine open neighbourhood $U=\Spec A$, so that $x$ corresponds to a prime ideal $\fp \subset A$. Because $X$ is integral, $A$ is cancellative, as is its localization $\cO_{X,x}=A_\fp$. Thus, $X$ is cancellative and, in particular, pc.
 
Suppose instead that $X$ is irreducible and pc. We first prove that $X$ is cancellative. Let $x \in X$ with affine open neighbourhood $U=\Spec A$. Since $X$ is irreducible, $A$ has a unique minimal ideal $\fp$ and since $U$ is pc, $A$ is pc by Lemma \ref{lem:pcProperties}\ref{part:pcLoc}. Then $A$ fails to be cancellative only if it has nontrivial zero divisor, say $ab =0$ for non-zero $a,b \in A$. In this case, the inverse image $\fp_b=\iota_b^{-1}(\fm)$ of the maximal ideal $\fm=A[b^{-1}]-A[b^{-1}]^\times$ of $A[b^{-1}]$ under the localization map $\iota_b:A\to A[b^{-1}]$ is a prime ideal that contains $a$ but not $b$. Similarly, there is a prime ideal $\fp_b$ that contains $b$ but not $a$. Since $ab=0$, the intersection $\fp_a\cap\fp_b$ cannot contain a prime ideal. This contradicts the fact that $A$ has a unique minimal ideal. Hence, $A$ and all of its localizations are cancellative. This proves that $X$ is cancellative.
 
To complete the proof that $X$ is integral, denote by $\eta$ the generic point of $X$. Lemma \ref{lem:genericSubmonoid} implies that $\cO_X(U)$ is a submonoid of the cancellative stalk $\cO_{X,\eta}$ for every open subset $U$ of $X$, and therefore is itself cancellative. Thus $X$ is integral.
\end{proof}

\begin{Def}
\begin{enumerate}[wide,labelwidth=!, labelindent=0pt,label=(\roman*)]
\item A valuation monoid is a cancellative pointed monoid $A$ such that $\Frac A=\{a,a^{-1}\mid a\in A\}$.

\item A morphism $\varphi:Y\to X$ of monoid schemes is proper if for all valuation monoids $A$ with inclusion $\iota:A\to\Frac A$ and all morphisms $\mu:\Spec\Frac A\to Y$ and $\nu:\Spec A\to X$ with $\varphi\circ\mu=\nu\circ\iota^\ast$, there exists a unique morphism $\hat\nu:\Spec A\to Y$ such that the diagram
\[
 \begin{tikzcd}[column sep=60pt,row sep=20pt]
  \Spec \Frac A \ar[r,"\mu"] \ar[d,swap,"\iota^\ast"]      & Y \ar[d,"\varphi"] \\
  \Spec A \ar[r,swap,"\nu"] \ar[dashed,ur,"\hat\nu"] & X
 \end{tikzcd}
\]
commutes.

\item A monoid scheme $X$ of finite type is proper if the terminal morphism $X\to\Spec\Fun$ is proper.
\end{enumerate}
\end{Def}

\begin{Prop}
\label{prop:intPropAut}
Let $X$ be a proper, integral and torsion free monoid scheme of finite type. Then $\Gamma X = \Fun$.
\end{Prop}

\begin{proof}
 Let $\eta$ be the generic point of $X$. Since $X$ is integral, all stalks $\cO_{X,x}$, $x \in X$, are submonoids of the generic stalk $\cO_{X,\eta}$; see Lemma \ref{lem:genericSubmonoid}. Thus
 \[
  \Gamma X \ = \ \bigcap_{x\in X} \cO_{X,x},
 \]
the intersection being taken in $\cO_{X,\eta}$. 
 
 In order to prove that $\Gamma X = \Fun$, we consider an element $f\in\cO_{X,\eta}$ and assume that $f\notin\{0,1\}$. Being a proper scheme, $X$ is of finite type, and since it is torsion free, we have $\cO_{X,\eta}\simeq \Fun[T_1^{\pm1},\dotsc,T_n^{\pm1}]$ for some $n\geq0$. Under this isomorphism, $f$ corresponds to a Laurent monomial $\prod_{i=1}^n T_i^{e_i}$ for some tuple $(e_1,\dotsc,e_n)\neq (0,\dotsc,0)$ in $\Z^n$. Therefore, we find a pointed monoid morphism
 \[
  v: \ \cO_{X,\eta}\simeq \Fun[T_1^{\pm1},\dotsc,T_n^{\pm1}] \ \longrightarrow \ \Fun[T^{\pm1}]
 \]
 that maps $f$ to $T^{-i}$ for some $i>0$. Since $X$ is of finite type, $U_\eta=\Spec\cO_{X,\eta}$ is an open subscheme of $X$.
 
 Let $\iota:\Fun[T]\to\Fun[T^{\pm1}]$ be the canonical inclusion and consider the diagram
 \[
  \begin{tikzcd}
   \Spec \Fun[T^{\pm1}] \ar{r}{v^\ast} \ar{d}[swap]{\iota^\ast} & U_\eta \ar[right hook->]{r} & X \ar{d}{\varphi} \\
   \Spec\Fun[T] \ar{rr}[below]{\nu} \ar[dashed]{urr}{\hat\nu}          &                             & \Spec\Fun
  \end{tikzcd}
 \]
 whose outer square commutes as $\Spec\Fun$ is terminal and where $\hat\nu$ is the unique morphism given by the defining property of the proper scheme $X$. Let $z=\gen T$ be the closed point of $\Spec\Fun[T]$ and $x=\hat\nu(z)$. Consider the induced morphism of stalks $\hat\nu_z^\#:\cO_{X,x}\to\cO_{\Spec\Fun[T],z}=\Fun[T]$. Since $v(f)=T^{-i}$ for $i>0$, the element $f\in\cO_{X,\eta}$ is not contained in $\cO_{X,x}$, which shows that $f$ is not a global section. This shows that $\Gamma X = \{0,1\}$, as desired.
\end{proof}

\subsection{Vector bundles}
\label{sec:vectF1Sch}

Let $X$ and $F$ be monoid schemes. A fibre bundle on $X$ with fibre $F$, or simply an $F$-bundle on $X$, is a morphism $\pi:E\to X$ of monoid schemes such that there is an open covering $\{U_i\}$ of $X$ and isomorphisms $\varphi_i:E\times_XU_i\to F \times U_i$, called trivializations, such that each diagram
\[
 \begin{tikzcd}
  E\times_XU_i \ar[rr,"\varphi_i","\sim"'] \ar[dr,swap,"\pi_i"] &     & F\times U_i \ar[dl,"\pr_{U_i}"]\\
                                                   & U_i 
 \end{tikzcd}
\]
commutes, where $\pi_i=E\times_XU_i\to U_i$ is the restriction of $\pi$ to $U_i$. Sometimes we suppress the morphism $\pi$ from the notation and say that $E$ is an $F$-bundle on $X$. An $F$-bundle $E$ on $X$ is trivializable 
if there existis a trivialization $E\simeq F\times X$. A morphism of $F$-bundles $E$ and $E'$ on $X$ is a commutative diagram
\[
 \begin{tikzcd}
  E \ar[rr] \ar[dr,"\pi"'] &   & E' \ar[dl,"\pi'"] \\
                           & X 
 \end{tikzcd}
\]
of morphisms of monoid schemes.

\begin{Rem}\label{rem:TrivialBundles}
 If $X=\Spec A$ is affine, then it has a unique closed point and thus every covering $\{U_i\}$ is trivial in the sense that $U_i=X$ for some $i$. Therefore every $F$-bundle on $X$ is trivializable.

As a consequence, given an $F$-bundle $\pi:E\to X$ on an arbitrary monoid scheme $X$, we can find trivializations $\varphi_i:E\times_XU_i\to F\times U_i$ for every chosen affine open covering $\{U_i\}$ of $X$. In particular, this holds for the minimal affine open covering $\{\Spec\cO_{X,x}\mid x\in X\text{ closed}\}$ if $X$ has enough closed points.
\end{Rem}

Let $n \in \mathbb{Z}_{\geq 0}$. A vector bundle on $X$ of rank $n$ is an $\A^n_\Fun$-bundle. We denote the category of finite rank vector bundles on $X$, together with all bundle morphisms, by $\Vect(X)$.

\begin{Rem} 
In contrast to vector bundles on schemes over a field, we need not require any additional datum to describe vector bundles, since the `vector space structure' of $\A^n_\Fun$ is intrinsically given, and coordinate changes of an $\A^n_\Fun$-bundle are necessarily `$\Fun$-linear'. This follows from the fact that every $A$-linear automorphism of $A[T_1,\dotsc,T_n]$ is graded; cf.\ the proof of Proposition \ref{prop:VectLocFree}.
\end{Rem}

As in algebraic geometry over a field, vector bundles correspond to finite locally free sheaves. We briefly review the definitions.

Let $X$ be a monoid scheme. An $\cO_X$-module is a sheaf $\cF$ of pointed sets on $X$ together with a morphism $\cO_X\times\cF\to\cF$ of sheaves such that $\cF(U)$ is an $\cO_X(U)$-module and the restriction maps $\cF(U)\to\cF(V)$ are pointed $\cO_X(U)$-module homomorphisms for all open subsets $V\subset U \subset X$. A morphism of $\cO_X$-modules is a morphism $\varphi:\cF\to\cF'$ of sheaves such that $\cF(U)\to\cF'(U)$ is a pointed $\cO_X(U)$-module homomorphism for every open $U \subset X$. This defines the category $\cO_X \Modu$ of $\cO_X$-modules on $X$.

An $\cO_X$-module $\cF$ is said to be finite locally free if every point $x \in X$ has an open neighbourhood $x \in U \subset X$ such that $\cF \vert_U$ is a free $\cO_X \vert_U$-module of finite rank. We denote by $\LF(X)$ the full subcategory of $\cO_X \Modu$ on finite locally free sheaves.

The relation between finite locally free sheaves and vector bundles uses the symmetric algebra. Let $A$ be a pointed monoid and $M$ an $A$-module. The symmetric algebra of $M$ is the $\mathbb{Z}_{\geq 0}$-graded pointed monoid
\[
 \Sym(M) \ = \ \bigoplus_{i\in \mathbb{Z}_{\geq 0}} \, \Sym_i(M) 
\]
 where
 \[
  \Sym_i(M) \ = \ M^{\otimes i} \, / \, \gen{a_1\otimes\dotsb\otimes a_i=a_{\sigma(1)}\otimes\dotsb\otimes a_{\sigma(i)} \mid a_1,\dotsc,a_i\in M, \sigma \in \Sigma_i}
 \]
 for $i>0$ and $\Sym_0(M)=A$. The multiplication of $\Sym(M)$ is given by concatenation of tensors and the inclusion $A=\Sym_0(M)\hookrightarrow\Sym(M)$ is a monoid morphism.

\begin{Prop}\label{prop:VectLocFree}
 The sheafification of the functor $\Spec\circ\Sym\circ\Gamma(X,-)$ defines an equivalence of categories $\vect:\LF(X)\to\Vect(X)$.
\end{Prop}

\begin{proof}
 This is proven similarly to the corresponding fact for schemes over a field. We briefly sketch the key arguments, but forgo to verify all details. First of all, note that $\vect(\cF)$ is indeed a vector bundle since 
 \[
  \Sym(A^{\oplus n}) \ \simeq \ A[T_1,\dotsc,T_n] \ = \ \Fun[T_1,\dotsc,T_n]\otimes_\Fun A
 \]
 and therefore $\vect(U)=\Spec\Sym\big(\Gamma(U,\cF)\big)\simeq \A^n_\Fun\times\Spec\Gamma U$ for every affine open $U$ of $X$. 
 
 A quasi-inverse $\lf:\Vect(X)\to\LF(X)$ of $\vect$ can be defined as follows. Let $\pi:E\to X$ be a vector bundle of rank $n$, $U=\Spec A$ an affine open of $X$ and $V=E\times_XU$. A trivialization $\varphi_U:V \xrightarrow[]{\sim} U\times_\Fun\A^n_\Fun$ defines an isomorphism
 \[
  f_U: \Gamma V \ \simeq \ A\otimes_\Fun\Fun[T_1,\dotsc,T_n] \ \simeq \ A[T_1,\dotsc,T_n].
 \]
 Since the $A$-linear automorphisms of $A[T_1,\dotsc,T_n]$ correspond to the images of $T_1,\dotsc,T_n$, which must be of the form $f_U(T_i)=a_iT_{\sigma(i)}$ for some permutation $\sigma\in S_n$ and some $a_1,\dotsc,a_n\in A^\times$, the $A$-invariant subsets 
 \[
  \Gamma V_{1,i} \ = \ f_U^{-1}\big(\{aT_i\mid a\in A\}\big) \qquad \text{and} \qquad \Gamma V_1 \ = \ \bigcup_{i=1}^n \Gamma V_{1,i}
 \]
of $\Gamma V$ do not depend on the choice of trivialization $\varphi_U$ up to a permutation of indices. This yields a canonical representation of $\Gamma V$ as a symmetric algebra $\Sym(\Gamma V_1)$.
 
Define $\lf(E)(U)$ to be the set of $A$-linear maps $s:\Gamma V_1\to A$ such that $s(\Gamma V_{1,i})=\{0\}$ for all but one $i\in\{1,\dotsc,n\}$,
 which is a free $A$-module. The sheafification of the assignment $U\mapsto \lf(E)(U)$ defines a functor $\lf:\Vect(X)\to\LF(X)$ that is quasi-inverse to $\vect$. 
\end{proof}

In light of Proposition \ref{prop:VectLocFree}, we allow ourselves to consider vector bundles as sheaves. Note that under the correspondence $\Vect(X)\to\LF(X)$, line bundles (vector bundles of rank one) correspond to invertible sheaves, that is, locally free sheaves of rank one.

Define $\Pic X$ to be the set of isomorphism classes of invertible sheaves on $X$ together with the group operation induced by $\otimes$. By abuse of language, we call elements of $\Pic X$ line bundles and sometimes identify an isomorphism class with a chosen representative. The neutral element of $\Pic X$ is the class of $\cO_X$ and the inverse of a line bundle $\cL$ is the dual line bundle $\cL^\vee=\HHom_{\cO_X}(\cL,\cO_X)$ where $\HHom_{\cO_X}(\cL,\cO_X)$ is the sheafification of the functor $U\mapsto\Hom_{\Gamma U}\big(\cL(U),\cO_X(U)\big)$.

\subsection{Locally projective sheaves}
\label{sec:LPsheaves}

Let $U=\Spec A$ be an affine monoid scheme. An $A$-module $M$ defines an $\cO_U$-module $\widetilde M$ with $\widetilde M(U_h)=M\otimes_A A[h^{-1}]$. For an arbitrary monoid scheme $X$, we say that an $\cO_X$-module $\cF$ is finite locally projective if there exists an open covering $\{U_i\}$ of $X$ such $M_i=\cF(U_i)$ is a finitely generated projective $\cO_X(U_i)$-module and such that $\cF\vert_{U_i}\simeq\widetilde{M_i}$ as sheaves on $U_i$.
We denote the category of finite locally projective sheaves by $\LP(X)$. Finite locally free sheaves are locally projective sheaves. The converse implication holds for the following class of monoid schemes.

\begin{Lem}
 \label{lem:LocProjLocFree}
 Let $X$ be a pc monoid scheme. Then every finite locally projective sheaf is finite locally free.
\end{Lem}

\begin{proof}
 Let $\cF$ be a finite locally projective sheaf on $X$ and $\{U_i\}$ an affine open covering such that $\cF(U_i)$ is a finitely generated projective $\cO_X(U_i)$-module. Then $U_i=\Spec\cO_{X,x_i}$ is pc where $x_i$ is the unique closed point of $U_i$. If $e^2=e=1\cdot e$ in a pc pointed monoid, then either $e=0$ or $e=1$. Thus, by Lemma \ref{lem:ProjMod}, $\cF(U_i)$ is a free $\cO_X(U_i)$-module of finite rank, which shows that $\cF$ is finite locally free.
\end{proof}

A morphism of $\cO_X$-modules $f: \mathcal{E} \rightarrow \mathcal{F}$ is called normal if the pointed $\mathcal{O}_{X,x}$-module homomorphism $f_x: \mathcal{E}_x \rightarrow \mathcal{F}_x$ is normal for each $x \in X$. The category $\cO_X \Modu$ has a proto-exact structure in which conflations are kernel-cokernel pairs of normal morphisms \cite[Lemma 5.6]{chu2012}; see also \cite[Proposition 3.13]{jun2020}. Normal morphisms define a proto-exact subcategory $\cO_X \Modu^{\nor}$ of $\cO_X \Modu$.

The category $\LP(X)$ is an extension closed subcategory of $\cO_X \Modu$ and so inherits a proto-exact structure. The full subcategory $\LP^\nor(X)$ of $\cO_X \Modu^\nor$ on finite locally projective sheaves also has aninduced proto-exact structure.

We write $\oplus$ and $\otimes$ 
for the direct sum and tensor product on $\cO_X \Modu$, respectively. Both $\oplus$ and $\otimes$ induce bifunctors on $\LP(X)$, making $\cO_X \Modu$ and $\LP(X)$ into symmetric bimonoidal categories, as well as the respective subcategories $\cO_X \Modu^{\nor}$ and $\LP^\nor(X)$.

\begin{Lem}[{\cite[Theorem 5.12]{chu2012}}]
\label{lem:VectUniqSplit}
Let $X$ be a monoid scheme. Then the proto-exact category $\LP(X)$ is uniquely split and combinatorial.
\end{Lem}

Let $\angleser{\mathcal{O}_X} \subset \LP^\nor(X)$ be the full subcategory of objects which are isomorphic to $\mathcal{O}_X^{\oplus n}$ for some $n \in \mathbb{Z}_{\geq 0}$. It is a proto-exact symmetric bimonoidal subcategory. We can form the proto-exact symmetric bimonoidal category 
\[
\angleser{\mathcal{O}_X}[\Pic(X)]
:=
\bigoplus_{\mathcal{M} \in \Pic(X)} \angleser{\mathcal{O}_X}.
\]
See Appendix \ref{app:DirectSumCategories} for the definition of the right hand side.

The following result plays an important role in the remainder of the paper.

\begin{Prop}
\label{prop:VectDecomp}
Let $X$ be an integral monoid scheme. Then the functor
\[
F:= \bigoplus_{\mathcal{M} \in \Pic(X)} (-)_{\mathcal{M}} \otimes\mathcal{M}:
\angleser{\mathcal{O}_X}[\Pic(X)]
\rightarrow
\LP(X)
\]
is symmetric bimonoidal, proto-exact, essentially surjective and bijective on inflations and deflations. In particular, every finite locally projective sheaf on $X$ decomposes uniquely into a direct sum of line bundles.
\end{Prop}

\begin{proof}
That $F$ is symmetric bimonoidal is a direct calculation. It is clear that $F$ is proto-exact. By Proposition \ref{prop:intVsCanc}, the scheme $X$ is irreducible and pc. Proposition \ref{lem:LocProjLocFree} therefore implies that finite locally projective sheaves on $X$ are finite locally free. Since $X$ is irreducible, it has a unique generic point. The proof of \cite[Theorem 5.14]{chu2012} then applies to show that any finite locally projective sheaf on $X$ is isomorphic to a direct sum of line bundles, from which it follows that $F$ is essentially surjective.

To see that $F$ is bijective on inflations and deflations, consider the map
\[
 \Hom\Big(\bigoplus_{k=1}^n \mathcal{E}_{\mathcal{M}_{i_k}}, \ \bigoplus_{l=1}^m \mathcal{F}_{\mathcal{M}_{j_l}}\Big) \ \longrightarrow \ \Hom\Big(\bigoplus_{k=1}^n \mathcal{E}_{\mathcal{M}_{i_k}} \otimes \mathcal{M}_{i_k}, \ \bigoplus_{l=1}^m \mathcal{F}_{\mathcal{M}_{j_l}} \otimes \mathcal{M}_{j_l}\Big)
\]
whose domain and codomain is a set of morphisms in $\angleser{\mathcal{O}_X}[\Pic(X)]$ and $\LP(X)$, respectively. The left hand side is
\[
\prod_{k,l,i_k=j_l}
\Hom_{\angleser{\mathcal{O}_X}}(\mathcal{E}_{\mathcal{M}_{i_k}},\mathcal{F}_{\mathcal{M}_{j_l}}).
\]
Since $\LP(X)$ is split (Lemma \ref{lem:VectUniqSplit}), the subset of inflations or deflations in the codomain is
\begin{multline*}
\Hom^{\textnormal{inf}/\textnormal{def}}_{\LP(X)}(\bigoplus_{k=1}^n \mathcal{E}_{\mathcal{M}_{i_k}} \otimes \mathcal{M}_{i_k},\bigoplus_{l=1}^m \mathcal{F}_{\mathcal{L}_{j_l}} \otimes \mathcal{M}_{j_l})
\simeq\\
\prod_{\substack{k,l \\ i_k=j_l}}
\Hom^{\textnormal{inf}/\textnormal{def}}_{\LP(X)}(\mathcal{E}_{\mathcal{M}_{i_k}} \otimes \mathcal{M}_{i_k},\mathcal{F}_{\mathcal{M}_{j_l}} \otimes \mathcal{M}_{j_l}).
\end{multline*}
It therefore suffices to consider the case $n=m=1$ and $i_1=j_1$. In this case, $
\Hom_{\angleser{\mathcal{O}_X}[\Pic(X)]}(\mathcal{E}_{\mathcal{M}},\mathcal{F}_{\mathcal{M}}) = \Hom_{\LP(X)}(\mathcal{E}_{\mathcal{M}}, \mathcal{F}_{\mathcal{M}})$, which we claim can be identified with $\Hom_{\LP(X)}(\mathcal{E}_{\mathcal{M}} \otimes \mathcal{M},\mathcal{F}_{\mathcal{M}} \otimes \mathcal{M})$. To see this, consider the quasi-inverse auto-equivalences
\[
- \otimes \mathcal{M}:
\cO_X \Modu
\leftrightarrows
\cO_X \Modu
: - \otimes \mathcal{M}^{\vee}.
\]
Let $f: \mathcal{E} \rightarrow \mathcal{F}$ be a normal morphism in $\cO_X \Modu$. The stalk morphism $(f \otimes \mathcal{M})_x$ is naturally identified with
\[
f_x \otimes \id_{\mathcal{M}_x}:
\mathcal{E}_x \otimes_{\mathcal{O}_{X,x}} \mathcal{M}_x
\rightarrow
\mathcal{F}_x \otimes_{\mathcal{O}_{X,x}} \mathcal{M}_x,
\]
which is normal. It follows from this, and the fact that the quasi-inverse of $-\otimes \mathcal{M}$ is of the same form, that $-\otimes \mathcal{M}$ restricts to an autoequivalence of $\LP(X)$. This proves the claim and finishes the proof of the proposition.
\end{proof}

\begin{Ex}
The assumption that $X$ is irreducible in Proposition \ref{prop:VectDecomp} cannot be dropped. Indeed, let $X=\Proj\big(\Fun[T_0,T_1,T_2]/\gen{T_0T_1T_2=0}\big)$, which is the union of the three coordinate lines in $\P^2_\Fun$. The three canonical opens of $\P^2_\Fun$ provide a covering $X=U_0\cup U_1\cup U_2$, where $U_i\simeq \Spec\big(\Fun[T_j,T_k]/\gen{T_jT_k=0}\big)$ consists of the three points $\gen{T_j}$, $\gen{T_k}$ and $\gen{T_j,T_k}$ if $\{i,j,k\}=\{0,1,2\}$ and where $U_i\cap U_j\simeq\Spec\big(\Fun[T_k^{\pm1}]\big)$ consists of a single point $\gen{T_k}$. The situation is illustrated in Figure \ref{fig:ExNondecBundle}.
\begin{figure}[htb]
 \centering
 \begin{tikzpicture}[x=1.0cm,y=1.0cm]
  \begin{scope}[blend mode=multiply]
   \fill[gray!20!white]  (90:0.8) circle (1.6); 
   \fill[gray!20!white] (210:0.8) circle (1.6); 
   \fill[gray!20!white] (330:0.8) circle (1.6); 
  \end{scope}
  \filldraw[gray!90!white,rotate around={0:(0,0)}=0] (270:1cm) ellipse (1.7 and 0.1); 
  \filldraw[gray!90!white,rotate around={120:(0,0)}=0] (270:1cm) ellipse (1.7 and 0.1); 
  \filldraw[gray!90!white,rotate around={240:(0,0)}=0] (270:1cm) ellipse (1.7 and 0.1); 
  \filldraw  (90:2cm) circle (2pt); 
  \filldraw (210:2cm) circle (2pt); 
  \filldraw (330:2cm) circle (2pt); 
  \node[font=\footnotesize] (T0) at (150:1.4) {$\gen{T_0}$};
  \node[font=\footnotesize] (T1) at  (30:1.4) {$\gen{T_1}$};
  \node[font=\footnotesize] (T2) at (270:1.4) {$\gen{T_2}$};
  \node[font=\footnotesize] (T01) at ( 70:2.1) {$\gen{T_0,T_1}$};
  \node[font=\footnotesize] (T02) at (203:2.6) {$\gen{T_0,T_2}$};
  \node[font=\footnotesize] (T12) at (337:2.6) {$\gen{T_1,T_2}$};
  \node (T01) at (115:1.9) {$U_2$};
  \node (T01) at (235:1.9) {$U_1$};
  \node (T01) at (305:1.9) {$U_0$};
 \end{tikzpicture}
 \caption{The canonical covering of $\Proj\big(\Fun[T_0,T_1,T_2]/\gen{T_0T_1T_2=0}\big)$}
 \label{fig:ExNondecBundle}
\end{figure}

Consider the locally projective sheaf $\cF$ on $X$ with $\cF(U_i)\simeq\big(\Gamma(\cO_X,U_i))^{\oplus2}$, $i=0,1,2$, and whose transition functions $\varphi_{i,j}:\cF\vert_{U_i}(U_i\cap U_j) \to \cF\vert_{U_j}(U_i\cap U_j)$ are given by the permutation matrices
\[
 \varphi_{0,1} \ = \ \varphi_{0,2} \ = \ \begin{pmatrix} 1 & 0 \\ 0 & 1 \end{pmatrix} \qquad  \text{and} \qquad \varphi_{1,2} \ = \ \begin{pmatrix} 0 & 1 \\ 1 & 0 \end{pmatrix}. 
\]
Then $\cF$ is not isomorphic to a direct sum of line bundles.
\end{Ex}

Next, we describe the proto-exact summands of $\angleser{\mathcal{O}_X}[\Pic(X)]$.

\begin{Prop}
\label{prop:lineSubCat}
Let $X$ be an integral monoid scheme. Then $\Gamma X$ is cancellative and the global sections functor
\[
\Gamma(X,-)=\Hom_{\LP^\nor(X)}(\mathcal{O}_X,-):\angleser{\mathcal{O}_X}\rightarrow \Gamma X\Modu
\]
induces an exact equivalence $\angleser{\mathcal{O}_X}\xrightarrow[]{\sim}\Gamma X\proj^{\nor}$.
\end{Prop}

\begin{proof}
The functor $\Gamma(X,-)$ is $\oplus$-monoidal and hence exact, since the domain and codomain categories are uniquely split. The essential image of $\Gamma(X,-)$ consists of all finitely generated free $\Gamma X$-modules. Since $X$ is integral, the pointed monoid $\Gamma X$ is cancellative, being a submonoid of the stalk $\cO_{X,\eta}$ of the generic point $\eta$ of $X$ (Lemma \ref{lem:genericSubmonoid}).
Hence, projective $\Gamma X$-modules are free. That $\Gamma(X,-): \angleser{\mathcal{O}_X}
\rightarrow
\Gamma X \proj^{\nor}$ is fully faithful follows from the observation that both $\Hom_{\Gamma X \modu^{\nor}}(\Gamma X^{\oplus n}, \Gamma X^{\oplus m})$ and $\Hom_{\angleser{\mathcal{O}_X}} (\mathcal{O}_X^{\oplus n}, \mathcal{O}_X^{\oplus m})$ can be identified with the set of $\Gamma X$-valued $m \times n$-matrices with at most one non-zero entry in each row and column.
\end{proof}

\subsection{$K$-theory of monoid schemes}
\label{sec:KThyF1Sch}

In this section we describe the algebraic $K$-theory space of integral monoid schemes.

Let $X$ be a monoid scheme. Denote by $\mathcal{K}(X)=\mathcal{K}(\LP(X))$ the algebraic $K$-theory space of $\LP(X)$ and set $K_i(X) = \pi_i \mathcal{K}(X)$.

\begin{Thm}
\label{thm:KThySchemes}
Let $X$ be an integral monoid scheme. Then there is a homotopy equivalence
\[
\mathcal{K}(X)
\simeq
\sideset{}{'}\prod_{\mathcal{M} \in \Pic(X)} \mathbb{Z} \times B (\Gamma X^{\times} \wr \Sigma_{\infty})^+.
\]
\end{Thm}

\begin{proof}
By Lemma \ref{lem:reflectExactIso}, the functor $F$ from Proposition \ref{prop:VectDecomp} defines a homotopy equivalence $\mathcal{K}(\angleser{\mathcal{O}_X} [\Pic(X)])
\xrightarrow[]{\sim}
\mathcal{K}(X)
$. Using Propositions \ref{prop:DirectSumIsFilteredColimit} and \ref{prop:RestrictedProductIsFilteredColimit} and the fact that $K$-theory commutes with filtered colimits and finite direct sums of categories (see \cite[\S 2]{quillen1973}), we obtain
\begin{multline*}
\mathcal{K}(X)
\simeq
\varinjlim_{S\in \mathscr{P}_{< \infty}(\Pic(X))}\mathcal{K}
\big(
\bigoplus_{s\in S}\angleser{\mathcal{O}_X}
\big)
\simeq
\varinjlim_{S\in \mathscr{P}_{< \infty}(\Pic(X))}\prod_{s\in S}\mathcal{K}(\angleser{\mathcal{O}_X})
\\
\simeq
\sideset{}{'}\prod_{\mathcal{M} \in \Pic(X)} \mathcal{K}(\angleser{\mathcal{O}_X}).
\end{multline*}
Using Proposition \ref{prop:lineSubCat}, we have $\mathcal{K}(\angleser{\mathcal{O}_X})
\simeq
\mathcal{K}(\Gamma X\proj^{\nor})$. The desired homotopy equivalence now follows from Theorem \ref{thm:KThyRPC}.
\end{proof}

\begin{Cor}
\label{cor:algKThyConnProp}
Let $X$ be a proper integral torsion free monoid scheme. Then there is a homotopy equivalence
\[
\mathcal{K}(X)
\simeq
\sideset{}{'}\prod_{\mathcal{M} \in \Pic(X)} \mathbb{Z} \times (B \Sigma_{\infty})^+.
\]
\end{Cor}

\begin{proof}
This follows from Propositions \ref{prop:intVsCanc} and \ref{prop:intPropAut} and Theorem \ref{thm:KThySchemes}.
\end{proof}

Next, we deduce some results at the level of $K$-theory groups. In particular, we recover \cite[Theorem 5.14]{chu2012}.

\begin{Cor}
\label{cor:K0Pic}
Le $X$ be an integral monoid scheme with $\Gamma X^{\times}$ finite. Then there is an isomorphism of graded rings
\[
K_{\bullet}(X) \simeq \pi_{\bullet}^s(B (\Gamma X)^{\times}_+)[\Pic(X)],
\]
where the right hand side is the group algebra of $\Pic(X)$ with coefficients in the graded ring  $\pi_{\bullet}^s(B (\Gamma X^{\times})_+)$.
\end{Cor}

\begin{proof}
At the level of abelian groups, the statement follows from Theorem \ref{thm:KThySchemes} by taking homotopy groups. For the ring structure, we use that the functor $F$ of Proposition \ref{prop:VectDecomp} is symmetric bimonoidal. Denote by $\angleser{\mathcal{O}_X}_{\mathcal{M}}$ the direct summand of $\angleser{\mathcal{O}_X}$ concentrated in degree $\mathcal{M} \in \Pic(X)$. Then $\otimes$ restricts to biexact functors 
\[
\otimes: \angleser{\mathcal{O}_X}_{\mathcal{M}_1}\times \angleser{\mathcal{O}_X}_{\mathcal{M}_2}\rightarrow \angleser{\mathcal{O}_X}_{\mathcal{M}_1 \mathcal{M}_2},
\qquad
\mathcal{M}_i \in \Pic(X).
\]
In this way we obtain a commutative diagram of pairings of $K$-theory spectra
\begin{center}
\begin{tikzcd}
\mathbf{K}(\angleser{\mathcal{O}_X}_{\mathcal{M}_1})\times \mathbf{K}(\angleser{\mathcal{O}_X}_{\mathcal{M}_2})\arrow[r] \arrow[d] & \mathbf{K}(\angleser{\mathcal{O}_X}_{\mathcal{M}_1 \mathcal{M}_2}) \arrow[d] \\
\mathbf{K}( \angleser{\mathcal{O}_X}[\Pic(X)]) \times \mathbf{K}( \angleser{\mathcal{O}_X}[\Pic(X)])  \arrow[r]           & \mathbf{K}( \angleser{\mathcal{O}_X}[\Pic(X)]).       
\end{tikzcd}
\end{center}
Compare with \cite[\S IV.6.6]{weibel2013}. The remaining statements follow.
\end{proof}

\begin{Ex}
Fix $n \geq 1$ and let $\mathbb{P}_{\Fun}^n$ be the $n$-dimensional projective space over $\Fun$. It is a monoid scheme whose base change $\mathbb{P}_{\Fun}^n \times_{\Spec(\Fun)} \Spec(\mathbb{Z})$ is the $n$-dimensional projective space $\mathbb{P}^n_{\mathbb{Z}}$ over the integers. There is an isomorphism of groups $\Pic(\mathbb{P}_{\Fun}^n) \simeq \mathbb{Z}$ which sends $\mathcal{O}_{\mathbb{P}_{\Fun}^n}(1)$ to $1 \in \mathbb{Z}$; see \cite[Theorem 2.6]{bothmer2011} or \cite[\S 5.4.3]{chu2012}. Corollary \ref{cor:K0Pic} then gives the isomorphism $K_0(\mathbb{P}_{\Fun}^n) \simeq \mathbb{Z}[\mathbb{Z}]$. In this way, we recover \cite[Corollary 5.15]{chu2012}. In particular, $K_0(\mathbb{P}_{\Fun}^n)$ is independent of $n$, in stark contrast to the case of projective space over a field.

Note, however, that, as explained in \cite[Theorem 5.17]{chu2012}, the relations for the $K$-theory of projective space over a field can be recovered by a comparison morphism between the $K$-theory and the $G$-theory of $\P^n_\Fun$. Since $G$-theory does not extend to the Grothendieck--Witt theory, we forgo pursuing this viewpoint in the present paper.
\end{Ex}
\begin{Rem}
For an arbitrary monoid scheme $X$, we do not expect $\mathcal{K}(X)$, as defined above using the proto-exact category of vector bundles, to be the best definition of the algebraic $K$-theory space of $X$. This is for reasons similar to the case of schemes over a field, where algebraic $K$-theory defined using vector bundles has poor cohomological properties without imposing mild assumptions on $X$. See, for example, \cite[Corollary 3.9]{thomason1990}. Since we do not study cohomological properties of $\mathcal{K}(X)$ in this paper, we ignore this issue.
\end{Rem}

\subsection{Projective bundles}
\label{sec:projBun}

Let $X$ be a monoid scheme. A projective bundle on $X$ is a $\P_{\mathbb{F}_1}^n$-bundle $\pi:E\to X$ for some $n\geq0$.

As in algebraic geometry over a field, there is a correspondence between classes of vector bundles, or finite locally free sheaves, and projective bundles on $X$. Let $\cE=\lf(E)$ be a locally free sheaf of rank $r$ on $X$. The sheafification of the functor $\Proj\circ\Sym\circ\Gamma(-,\cE)$ defines a $\P_{\mathbb{F}_1}^{r-1}$-bundle $\pi:\P\cE\to X$; see \cite[\S 7]{cortinas2015} for details on the $\Proj$ construction. Its restriction to an affine open $U$ of $X$ is isomorphic to $\P^{r-1}_U$ such that $\pi_U$ commutes with the structure map $\P^{r-1}_U\to U$.

\begin{Rem}
 We note without proof that every projective bundle is of the form $\P\cE$ for some finite locally free sheaf $\cE$ on $X$. Two finite locally free sheaves $\cE$ and $\cE'$ define isomorphic projective bundles $\pi:\P\cE\to X$ and $\pi':\P\cE'\to X$ if and only if there is a locally free sheaf $\cL$ of rank one such that $\cE'\simeq\cE\otimes\cL$.
\end{Rem}

\begin{Lem}\label{lem:sectionsPn}
 Let $A$ be a pointed monoid and $n\geq0$. Then $\Gamma\P^n_A=A$.
\end{Lem}

\begin{proof}
 Consider the canonical affine open $U_i=\Spec A[\frac{T_j}{T_i}|j=0,\dotsc,n]$ of $\P^n_A=\Proj A[T_0,\dotsc,T_n]$ and $U_\eta=\bigcap_{i=0}^n U_i=\Spec A[\frac{T_j}{T_i}|i,j=0,\dotsc,n]$. Since $\P^n_A$ is covered by the $U_i$ and the restriction maps $\Gamma U_i\to \Gamma U_\eta$, $i=0,\dotsc,n$, are injective, we have $\Gamma\P^n_A=\bigcap_{i=0}^n \Gamma U_i$ as an intersection inside $\Gamma U_\eta$.

 In particular, a global section $f\in\Gamma\P^n_A$ is contained in $\Gamma U_\eta$ and is therefore of the form $a \prod_{i=0}^n T_i^{e_i}$ for some $a\in A$ and $(e_0,\dotsc,e_n)\in\Z^{n+1}$ with $\sum_{i=0}^n e_i=0$. Since $f\in\Gamma U_i$, we have $e_i\leq 0$ for all $i=0,\dotsc,n$. Thus $e_0=\dotsb=e_n=0$, which shows that $f\in A$ as claimed.
\end{proof}

\begin{Prop}
\label{prop:globFunProjBun}
 Let $X$ be an integral monoid scheme and $\pi:\P\cE\to X$ a projective bundle. Then the map $\pi^* : \Gamma X\rightarrow \Gamma\P\cE$ is an isomorphism of pointed monoids. 
\end{Prop}

\begin{proof}
 Since $X$ is integral,
 it has a unique generic point $\eta$, all restriction maps are injective and $\Gamma X$ is the intersection $\bigcap\Gamma U$ inside $\cO_{X,\eta}$, where $U$ varies over all affine open subschemes of $X$.
 
For every affine open $U \subset X$, the bundle $\P\cE\times_XU\simeq\P_U^n$ trivializes. Since $\eta$ is contained in every open subset of $X$, we conclude that $\P\cE$ is irreducible with unique generic point $\hat\eta$ mapping to $\eta$. Moreover, since every affine open subscheme $U$ of $X$ is integral, $\P\cE$ is covered by the integral subschemes $\P\cE\times_XU\simeq\P_U^n$, which shows that $\P\cE$ is integral. We conclude, using Lemma \ref{lem:genericSubmonoid}, that $\Gamma\P\cE$ embeds into $\cO_{\P\cE,\hat\eta}$ and equals the intersection $\bigcap_U \Gamma U$ where $U$ varies over an affine open covering of $X$.
 
To show that $\pi^*$ is injective,
suppose that $\pi^*(s)=\pi^*(t)$ for $s,t\in\Gamma X$. Let $V\simeq\Spec A$ be an affine open of $X$ and $U_{V,i}\simeq\Spec A[\frac{T_j}{T_i}|j=0,\dotsc,n]$ the canonical open of $\P\cE\times_XV\simeq\P_V^n$. Since $(\pi\vert_V)^*:A\to A[\frac{T_j}{T_i}|j=0,\dotsc,n]$ is injective and $\pi^*(s)\vert_{U_{V,i}}=\pi^*(t)\vert_{U_{V,i}}$, we conclude that $s\vert_V=t\vert_V$. Covering $X$ with affine opens $V$ yields the equality $s=t$.
 
We turn to surjectivity of $\pi^*$. Let $s\in\Gamma\P\cE$. By Lemma \ref{lem:sectionsPn}, the restriction $s\vert_U$ of $s$ to $U=\P\times_XV\simeq\P_A^n$ comes from a global section $t_V\in\Gamma V$ for every affine open $V$ of $X$. Since $(\pi\vert_V)^*:\Gamma V\to\Gamma U$ is injective, $t_V$ is unique, and therefore the collection $\{t_V\}$, where $V$ varies through all affine opens of $X$, glues to a unique global section $t$ of $\P\cE$ with $\pi^*(t)=s$. This shows that $\pi^*$ is surjective and completes the proof.
\end{proof}

\begin{Lem}
 \label{lem:SectionPE}
 Let $X$ be a monoid scheme and $\pi:\P\cE\to X$ a projective bundle. Then there is a canonical section $\sigma:X\to\P\cE$ such that $\sigma(x)$ is the generic point of the fibre $\pi^{-1}(x)$ for every $x\in X$.
\end{Lem}

\begin{proof}
 Choose an affine open covering $\{U_i=\Spec A_i\}$ of $X$ and define $V_i=U_i\times_X\P\cE$. Since $\P\cE$ trivializes over affine opens, $V_i$ is isomorphic to $\P^n_{A_i}=\Proj A_i[T_0,\dotsc,T_n]$, with $n$ the fibre dimension of $\P\cE$. For each $i$, the graded $A_i$-linear pointed monoid morphism
 \[
  \begin{array}{cccc}
   s_i: & A_i[T_0,\dotsc,T_n] & \longrightarrow & A_i[\widehat T] \\
        & T_j                 & \longmapsto     & \widehat T
  \end{array}
 \]
 defines a morphism
 \[
  \sigma_i=s_i^\ast: U_i \ = \ \P^0_{U_i} \ = \ \Proj A_i[\widehat T] \ \longrightarrow \ \Proj A_i[T_0,\dotsc,T_n] \ \simeq \ V_i.
 \]
Since the $s_i$ are invariant under $A_i$-linear automorphisms of $A_i[T_0,\dotsc,T_n]$, they do not depend on the choice of identification $V_i\simeq \P^n_{A_i}$, and therefore coincide on the intersections of the $U_i$ and glue to a canonical morphism $\sigma:X\to \P\cE$. Since $\pi^\#(U_i)\circ s_i$ is the identity on $A_i[\widehat T]$ for every $i$, the composition $\pi\circ\sigma$ is the identity on $X$, which shows that $\sigma$ is a section to $\pi$. The restriction of $\sigma$ to $\sigma_x:\{x\}\to\pi^{-1}(x)$ corresponds to the graded $k(x)$-monoid morphism
 \[
  \begin{array}{cccc}
   s_x: & k(x)[T_0,\dotsc,T_n] & \longrightarrow & k(x)[\widehat T] \\
        & T_j                  & \longmapsto     & \widehat T
  \end{array}
 \]
 where $k(x)=\cO_{X,x}/\mathfrak{m}_x$ is the ``residue field'' at $x$. Since $s_x^{-1}(\gen 0)=\gen 0$, we conclude that $\sigma_x(x)$ is the generic point in $\pi^{-1}(x)$. This completes the proof.
\end{proof}

\begin{Lem}
 \label{lem:PicPn}
 Let $A$ be a pointed monoid and $n\geq1$. Then $\Pic(\P^n_A)=\{\cO_{\P^n_A}(m)\mid m\in\Z\}$.
\end{Lem}

\begin{proof}
 As explained in the example at the end of Section \ref{sec:KThyF1Sch}, this result is known for $A=\Fun$. 
 The inclusion $i:\Fun\to A$ induces a morphism $\pi:\P^n_A\to\P^n_\Fun$ and a group homomorphism $\pi^\ast:\Pic(\P^n_\Fun)\to\Pic(\P^n_A)$. 
 
 Choose any pointed monoid morphism $p:A\to \Fun$, such as sending all units to $1$ and all other elements to $0$. Then $i\circ p=\id_\Fun$ and $p$ induces a section $\sigma:\P^n_\Fun\to\P^n_A$ of $\pi$ and a retract $\sigma^\ast:\Pic(\P^n_A)\to\Pic(\P^n_\Fun)$ of $\pi^\ast$. We conclude that $\pi^\ast$ is injective.
 
 We turn to the surjectivity of $\pi^\ast$. Let $U_i=\Spec\Fun[T_j/T_i]_{j=0,\dotsc,n}$ be the canonical open subsets of $\P^n_\Fun=\Proj \Fun[T_0,\dotsc,T_n]$ and 
 \[
  U_\eta \ = \ U_0\cap\dotsc\cap U_n \ = \ \Spec\Fun[T_j/T_i]_{i,j=0\dotsc,n}. 
 \]
Let $V_i=U_i\times_{\P^n_\Fun}\P^n_A$ and $V_\eta=U_\eta\times_{\P^n_\Fun}\P^n_A$. Then $\{V_i\}$ is an affine open covering of $\P^n$ and every line bundle $\cL$ on $P^n_A$ trivializes over this covering. This means that we get, for every $i=0,\dotsc,n$, a commutative diagram
 \[
  \begin{tikzcd}[column sep=60pt,row sep=15pt]
   \cL(V_i) \ar{r}{\varphi_i}[swap]{\sim} \ar[d,swap,"\res_{V_i,V_\eta}"] & \Gamma V_i \ar[d,"\iota_i"] \\
   \cL(V_\eta) \ar{r}{\sim}[swap]{\varphi_{i,\eta}}                       & \Gamma V_\eta
  \end{tikzcd}
 \]
 of $\Gamma V_i$-linear maps whose right vertical arrow is the canonical inclusion that comes from the inclusion $\Gamma U_i\to\Gamma U_\eta$. Since the $\Gamma V_i$-linear map $\varphi_{i,\eta}\circ\res_{V_i,V_\eta}$ is determined by the image of $1$, which is of the form $a\prod_{i=0}^n T_i^{e_i}$ for some $a\in A$ and $(e_0,\dotsc,e_n)\in\Z^{n+1}$ with $\sum_{i=0}^n e_i=0$, and $a$ is invertible (since $\iota_i=\varphi_{i,\eta}\circ\res_{V_i,V_\eta}\circ\varphi_i^{-1}$ is a localization of pointed monoids), we can assume that 
 \[
  \varphi_{i,\eta}\circ\res_{V_i,V_\eta}(1) \ = \ \prod_{i=0}^n T_i^{e_i},
 \]
 after replacing $\varphi_i$ by $a^{-1}\varphi^{-1}$. This shows that we can choose the trivalizations $\varphi_i$ so that they restrict to bijections $\Gamma U_i\to \Fun[T_j/T_i]_{j=0,\dotsc,n}$. This yields a line bundle $\cL'$ on $\P^n_\Fun$ with $\pi^\ast(\cL')=\cL$. Therefore, $\pi^\ast$ is surjective.
\end{proof}

The following result is a monoid-theoretic analogue of a well-known result for schemes over algebraically closed fields.

\begin{Thm}
 \label{thm:ProjLineBundles}
 Let $X$ be an irreducible monoid scheme and $\pi:\P\mathcal{E}\to X$ a projective bundle. 
 Then the map
 \[
  \begin{array}{cccc}
   \varphi: & \Pic(X) \times \mathbb{Z} & \longrightarrow & \Pic(\mathbb{P}\mathcal{E}) \\
            & (\mathcal{L},m)           & \longmapsto     & \pi^* \mathcal{L} \otimes \mathcal{O}_{\mathbb{P}\mathcal{E}}(m)
  \end{array}            
 \]
 is an isomorphism of abelian groups.
\end{Thm}

\begin{proof}
 Let $\eta$ be the generic point of $X$ and $U_\eta=\Spec\cO_{X,\eta}$, which comes with a canonical morphism $\iota_\eta:U_\eta\to X$. Define $V_\eta=U_\eta\times_X\P\cE$, which comes with the cartesian diagram
 \[
  \begin{tikzcd}[column sep=60pt,row sep=15pt]
   V_\eta \ar[r,"\iota"] \ar[d,swap,"\pi_\eta"] & \P\cE \ar[d,swap,"\pi"] \\
   U_\eta \ar[r,"\iota_\eta" swap]                   & X \ar[bend right=30pt]{u}[right]{\sigma}
  \end{tikzcd}
 \]
 where $\sigma:X\to \P\cE$ is the canonical section of $\pi$ from Lemma \ref{lem:SectionPE}. This yields a commutative diagram of group homomorphisms
 \[
  \begin{tikzcd}[column sep=60pt,row sep=15pt]
   \Pic V_\eta \ar[<-,r,"\iota^\ast"] \ar[<-,d,swap,"\pi_\eta^\ast"] & \Pic \P\cE \ar[<-,d,swap,"\pi^\ast"] \\
   \Pic U_\eta \ar[<-,r,"\iota_\eta^\ast" swap]                           & \Pic X \ar[<-,bend right=30pt]{u}[right]{\sigma^\ast}
  \end{tikzcd}  ,
 \]
 where $\Pic U_\eta$ is trivial since $U_\eta$ is affine. Thus $\iota^\ast\circ\pi^\ast=\pi_\eta^\ast\circ\iota_\eta^\ast=0$. Since $\sigma^\ast\circ\pi^\ast$ is the identity on $\Pic X$, we conclude that $\pi^\ast$ is injective.
 
 By Lemma \ref{lem:PicPn}, $\Pic V_\eta=\{\cO_{V_\eta}(n)\mid n\in\Z\}$. Thus, the assignment $\cO_{V_\eta}(n)\mapsto \cO_{\P\cE}(n)$ defines a group homomorphism $r:\Pic V_\eta\to\Pic\P\cE$, which is a section of $\iota^\ast$. This shows that $\iota^\ast$ is a surjection.
 
Consider a line bundle $\cL\in\Pic\P\cE$ in the kernel of $\iota^\ast$, so that $\iota^\ast(\cL)\simeq\cO_{V_\eta}$. Choose an affine open covering $\{U_i\}$ of $X$, so that $U_i=\Spec A_i$ for $A_i=\Gamma U_i$, and define $V_i=U_i\times_X\P\cE$. This defines an open covering $\{V_i\}$ of $\P\cE$. Since the $U_i$ are affine, $\P\cE\vert_{U_i}$ is trivializable, so that $\P\cE\vert_{U_i}\simeq\P^n_{A_i}$ where $n$ is the fibre dimension of $\pi$. By Lemma \ref{lem:PicPn}, we conclude that $\cL\vert_{V_i}\simeq\cO_{V_i}(m)$ for some $m\in\Z$. Since $\cO_{V_\eta}(m)\simeq \cL\vert_{V\eta}\simeq\cO_{V_\eta}$, we have $m=0$. We therefore obtain a commutative diagram
 \[
  \begin{tikzcd}[column sep=40pt,row sep=20pt]
   \cL(V_i) \ar[r,"\sim"] \ar[d,swap,"\res_{V_i,V_\eta}"] & \cO_{\P\cE}(V_i) \ar[r,"\sim"] \ar[d,"\res_{V_i,V_\eta}"] & \pi^\ast(\cO_{X})(V_i) \ar[r,"\sim"] \ar[d,"\res_{V_i,V_\eta}"] & \pi^\ast\big(\sigma^\ast(\cL)\big)(V_i) \ar[d,"\res_{V_i,V_\eta}"] \\
   \cL(V_\eta) \ar[r,"\sim"]                              & \cO_{\P\cE}(V_\eta) \ar[r,"\sim"]                         & \pi^\ast(\cO_{X})(V_\eta) \ar[r,"\sim"]                         & \pi^\ast\big(\sigma^\ast(\cL)\big)(V_\eta)
  \end{tikzcd}  
 \]
 for every $i$. Since the $V_i$ cover $\P\cE$, we conclude that $\cL=\pi^\ast\big(\sigma^\ast(\cL)\big)$ is in the image of $\pi^\ast$.
 
 Altogether, this shows that there is a canonically split short exact sequence
 \[
  \begin{tikzcd}[column sep=50pt,row sep=20pt]
   0 \ar[r] & \Pic X \ar[r,"\pi^\ast" swap] & \Pic\P\cE \ar[r,"\iota^\ast" swap] \ar[l,bend right=30pt,swap,"\sigma^\ast"] & \Pic V_\eta \ar[r] \ar[l,swap,bend right=30pt,"r"] & 0,
  \end{tikzcd}
 \]
 which induces the isomorphism
 \[
  \begin{array}{ccccc}
   \Pic X\times\Z & \stackrel\sim\longrightarrow & \Pic X\times\Pic V_\eta     & \stackrel{(\pi^\ast,r)}\longrightarrow & \Pic\P\cE \\
   (\cL,m)          & \longmapsto                  & \big(\cL,\cO_{V_\eta}(m)\big) & \longmapsto                            & \cL\otimes\cO_{\P\cE}(m)
  \end{array}
 \]
 of the claim of the theorem.
\end{proof}

\subsection{A projective bundle formula}
\label{sec:projBundK}

We combine our earlier results to prove a projective bundle formula for the $K$-theory space of a monoid scheme. This gives a monoid-theoretic analogue of Quillen's projective bundle formula over fields \cite[Theorem 2.1 of Section 8]{quillen1973}.

\begin{Thm}
\label{thm:projBunSpaces}
Let $X$ be an integral monoid scheme and $\pi:\P\mathcal{E}\to X$ a projective bundle. Then there is a homotopy equivalence
\[
\mathcal{K}(\mathbb{P} \mathcal{E})
\xrightarrow[]{\sim}
\sideset{}{'}\prod_{n \in \mathbb{Z}}\mathcal{K}(X).
\]
\end{Thm}

\begin{proof}
Combining Propositions \ref{prop:VectDecomp} and \ref{prop:globFunProjBun} and Theorem \ref{thm:ProjLineBundles}, we obtain a diagram of functors
\begin{multline*}
\LP^\nor(\mathbb{P} \mathcal{E})
\leftarrow
\angleser{\mathcal{O}_{\mathbb{P}\mathcal{E}}}[\Pic(\mathbb{P}\mathcal{E})]
\simeq
\angleser{\mathcal{O}_{X}}[\Pic(X)\times \mathbb{Z}]
\\
\simeq
(\angleser{\mathcal{O}_{X}}[\Pic(X)])[\mathbb{Z}]
\rightarrow
\LP(X)[\mathbb{Z}].
\end{multline*}
By Lemma \ref{lem:reflectExactIso}, the first and last functors induce homotopy equivalences of $K$-theory spaces. We therefore have homotopy equivalences
\[
\mathcal{K}(\mathbb{P} \mathcal{E})\simeq \mathcal{K}(\LP(X)[\mathbb{Z}])\simeq \sideset{}{'}\prod_{n \in \mathbb{Z}}\mathcal{K}(X).
\qedhere
\]
\end{proof}

\begin{Cor}
Let $X$ be an integral monoid scheme. Then there is an isomorphism of graded rings
\[
K_{\bullet}(\mathbb{P}\mathcal{E})
\simeq
K_{\bullet}(X) \otimes_{\mathbb{Z}} \mathbb{Z}[\mathbb{Z}].
\]
\end{Cor}



\section{Grothendieck--Witt theory of monoid schemes}
\label{sec:GWThySch}

\subsection{Duality for locally free sheaves}
\label{sec:vbDuality}
Let $X$ be a monoid scheme. There is a bifunctor
\[
\HHom_{\mathcal{O}_X}(-,-) : \cO_X \Modu^{\op} \times \cO_X \Modu \rightarrow \cO_X \Modu
\]
defined so that, for $\mathcal{E}, \mathcal{F} \in \cO_X \Modu$ and an open set $U \subset X$, we have
\[
\HHom_{\mathcal{O}_X}(\mathcal{E},\mathcal{F})(U) = \Hom_{\mathcal{O}_{X \vert U}}(\mathcal{E}_{\vert U}, \mathcal{F}_{\vert U}).
\]
Let $\mathcal{L}$ be a line bundle on $X$. As in the local case (Section \ref{sec:monoids}), and unlike the case of schemes over a field, the functor $\HHom_{\mathcal{O}_X}(-,\mathcal{L})$ does not define a duality functor on $\LF^{\nor}(X)$. This can be remedied as follows. Let $X$ be a pc monoid scheme. Given $\mathcal{E} \in \cO_X \Modu$, define a sheaf of pointed sets $P^{\mathcal{L}}(\mathcal{E})$ on $X$ by
\[
P^{\mathcal{L}}(\mathcal{E})(U)
=
\Hom^{\nor}_{\mathcal{O}_{X \vert U}}(\mathcal{E}_{\vert U},\mathcal{L}_{\vert U}),
\]
the right hand side being the set of normal $\mathcal{O}_{X \vert U}$-module homomorphisms $\mathcal{E}_{\vert U} \rightarrow \mathcal{L}_{\vert U}$. Because $X$ is pc, $P^{\mathcal{L}}(\mathcal{E})(U)$ has a natural $\mathcal{O}_X(U)$-module structure, as follows from a local calculation using Lemma \ref{lem:dualFreeAMod}.

\begin{Prop}
\label{prop:vectDuality}
Let $\mathcal{L}$ be a line bundle on a reversible pc monoid scheme $X$. Let $\Theta^{\mathcal{L}}: \id_{\LF^{\nor}(X)} \Rightarrow P^{\mathcal{L}} \circ (P^{\mathcal{L}})^{\op}$ be the natural isomorphism with components
\[
\Theta^{\mathcal{L}}_{\mathcal{E}}(s)(f) = f(s),
\qquad
s \in \mathcal{E}(U), \; f \in P^{\mathcal{L}}(\mathcal{E})(U)
\]
where $U \subset X$ is an open subset. Then $(\LF^{\nor}(X), P^{\mathcal{L}}, \Theta^{\mathcal{L}})$ is a uniquely split combinatorial proto-exact category with duality. Moreover, the Reduction Assumption (see Section \ref{sec:catWD}) holds.
\end{Prop}

\begin{proof}
That $P^{\mathcal{L}}$ sends locally free sheaves to locally free sheaves follows from the fact that $X$ is reversible and a local calculation using Lemma \ref{lem:dualFreeAMod}. Let $f: \mathcal{E} \rightarrow \mathcal{F}$ be a morphism in $\LF^{\nor}(X)$. For each $x \in X$, the stalk morphism $P^{\mathcal{L}}(f)_x : P^{\mathcal{L}}(\mathcal{F})_x \rightarrow P^{\mathcal{L}}(\mathcal{E})_x$ can be identified with
\[
(-) \circ f_x
:
\Hom_{\mathcal{O}_{X,x} \Modu^{\nor}}(\mathcal{F}_x,\mathcal{L}_x) \rightarrow \Hom_{\mathcal{O}_{X,x} \Modu^{\nor}}(\mathcal{E}_x,\mathcal{L}_x).
\]
Since the composition of normal morphisms is normal, $P^{\mathcal{L}}(f)$ is well-defined. That $P^{\mathcal{L}}$ is compatible with $\oplus$ follows from Lemma \ref{lem:dualFreeAMod}. That $\Theta^{\mathcal{L}}$ is a natural isomorphism follows from a local calculation using Proposition \ref{prop:projDuality}. 

It remains to verify the Reduction Assumption. Whether or not the induced map $\psi_{N \git U} : N \git U \rightarrow P(N \git U)$ is an isomorphism can be checked locally, in which case it reduces to Proposition \ref{prop:projDuality}.
\end{proof}

\begin{Rem}
Proposition \ref{prop:vectDuality} admits the following generalization, which can be seen as the natural commutative globalization of the setting of Proposition \ref{prop:projDuality}. Let $\sigma: X \rightarrow X$ be an involution and $\epsilon \in \Gamma(X,\mathcal{O}_X)^{\times}$. Assume that $\sigma^* \mathcal{L} \simeq \mathcal{L}$ and $\epsilon \sigma(\epsilon)=1$. Let $P^{\mathcal{L},\sigma} : \LF^{\nor}(X)^{\op} \rightarrow \LF^{\nor}(X)$ be the functor $\mathcal{E} \mapsto \HHom_{\mathcal{O}_X}^{\nor}(\sigma^* \mathcal{E}, \mathcal{L})$ and define $\Theta^{\mathcal{L},\sigma,\epsilon}$ by $\Theta^{\mathcal{L},\sigma,\epsilon}_{\mathcal{E}}(s)(f) = \epsilon \sigma(f(s))$. Then the analogue of Proposition \ref{prop:vectDuality} holds for $(\LF^{\nor}(X), P^{\mathcal{L},\sigma}, \Theta^{\mathcal{L},\sigma,\epsilon})$. The results which follow hold also at this level of generality, with essentially the same proofs. We note only that the involution of the pointed monoid $\Gamma X$ is determined through the isomorphism $\Gamma X \simeq \End_{\mathcal{O}_X \Modu}(\mathcal{M})$ via the formula $f \mapsto \psi_{\mathcal{M}} \circ P^{\mathcal{L},\sigma}(f) \circ \psi_{\mathcal{M}}^{-1}$. However, for ease of exposition, we restrict to the case $\sigma = \id_X$ and $\epsilon = 1$.
\end{Rem}

When $\mathcal{L} = \mathcal{O}_{X}$ we omit it from the notation so that, for example, $P^{\mathcal{O}_X}=P$.

\begin{Lem}
\label{lem:twistedParity}
Let $X$ be a reversible pc monoid scheme. If line bundles $\mathcal{L}$, $\mathcal{L}^{\prime}$ are equal in $\Pic(X) \slash \Pic(X)^2$, then there is an equivalence of proto-exact categories with duality
\[
(\LF^{\nor}(X), P^{\mathcal{L}}, \Theta^{\mathcal{L}}) \simeq (\LF^{\nor}(X), P^{\mathcal{L}^{\prime}}, \Theta^{\mathcal{L}^{\prime}}).
\]
\end{Lem}

\begin{proof}
Under the assumption of the lemma, there exists a line bundle $\tilde{\mathcal{L}} \in \LF^{\nor}(X)$ and an isomorphism
\begin{equation}
\label{eq:squareIso}
\mathcal{L} \otimes \tilde{\mathcal{L}} \xrightarrow[]{\sim} \tilde{\mathcal{L}}^{\vee} \otimes \mathcal{L}^{\prime}.
\end{equation}
Let $T$ be the exact autoequivalence $- \otimes \tilde{\mathcal{L}} : \LF^{\nor}(X) \rightarrow \LF^{\nor}(X)$; see the proof of Proposition \ref{prop:VectDecomp}. Then $(T,\eta)$ is an equivalence of categories with duality, where $\eta: T \circ P^{\mathcal{L}} \Rightarrow P^{\mathcal{L}^{\prime}} \circ T^{\op}$ is the natural isomorphism with components
\[
\eta_{\mathcal{E}}: P(\mathcal{E}) \otimes \mathcal{L} \otimes \tilde{\mathcal{L}} \rightarrow P(\mathcal{E}) \otimes \tilde{\mathcal{L}}^{\vee} \otimes \mathcal{L}^{\prime},
\qquad
\mathcal{E} \in \LF^{\nor}(X)
\]
defined using the chosen isomorphism \eqref{eq:squareIso}.
\end{proof}

\begin{Lem}
\label{lem:dualVsNorDual}
Let $X$ be a reversible pc monoid scheme and $\mathcal{E} \in \LF^{\nor}(X)$ a line bundle. Then there is a canonical isomorphism $\HHom_{\mathcal{O}_X}(\mathcal{E}, \mathcal{O}_X) \simeq P(\mathcal{E})$.
\end{Lem}

\begin{proof}
This follows from the local observation that, for a pc pointed monoid $A$, any $A$-module homomorphism $A \rightarrow A$ is normal.
\end{proof}

The functor $P$ is $\otimes$-monoidal, while $P^{\mathcal{L}}$ is not in general. Instead, $P^{\mathcal{L}}$ is $P$-monoidal, that is, there are coherent isomorphisms
\[
P^{\mathcal{L}}(\mathcal{E} \otimes \mathcal{F}) \simeq P(\mathcal{E}) \otimes P^{\mathcal{L}}(\mathcal{F}),
\qquad
\mathcal{E}, \mathcal{F} \in \LF^{\nor}(X)
\]
which are natural in $\mathcal{E}$ and $\mathcal{F}$. The functor $P^{\mathcal{L}}$ induces an involution
\begin{equation}
\label{eq:PicardInv}
P^{\mathcal{L}}: \Pic(X) \rightarrow \Pic(X).
\end{equation}
By Lemma \ref{lem:dualVsNorDual}, this map agrees with that induced by $\HHom_{\mathcal{O}_X}(-,\mathcal{L})$. We emphasize that, since $P^{\mathcal{L}}$ is only $P$-monoidal, the map \eqref{eq:PicardInv} is not a group homomorphism unless $\mathcal{L} \in \Pic(X)$ is trivial. Denote by $\Pic(X)^{P^{\mathcal{L}}}$ the $\mathbb{Z} \slash 2$-invariants of $\Pic(X)$.
The complement $\Pic(X)^* = \Pic(X) \setminus \Pic(X)^{P^{\mathcal{L}}}$ has a free $\mathbb{Z} \slash 2$-action.

Next, we study the compatibility of Proposition \ref{prop:VectDecomp} with duality. The subcategory $\angleser{\mathcal{O}_X} \subset \LF^{\nor}(X)$ is $P$-stable and so inherits from $\LF^{\nor}(X)$ a proto-exact duality, again denoted by $(P,\Theta)$. Define a proto-exact duality $(\mathbf{P}^{\mathcal{L}},\bm{\Theta}^{\mathcal{L}})$ on $\angleser{\mathcal{O}_X} [\Pic(X)]$ as follows. The functor $\mathbf{P}^{\mathcal{L}}$ is defined on basic objects by
\[
\mathbf{P}^{\mathcal{L}}(\mathcal{E}_{\mathcal{M}}) = P(\mathcal{E})_{P^{\mathcal{L}}(\mathcal{M})},
\qquad
\mathcal{E} \in \angleser{\mathcal{O}_X}, \; \mathcal{M} \in \Pic(X)
\]
and extended to $\angleser{\mathcal{O}_X} [\Pic(X)]$ by additivity. The natural isomorphism $\bm{\Theta}^{\mathcal{L}}$ has components $\bm{\Theta}^{\mathcal{L}}_{\mathcal{E}_{\mathcal{M}}} = \Theta^{\mathcal{L}}_{\mathcal{E}}$.

Note that, by Remark \ref{rem:CancelMonoidScheme} and Proposition \ref{prop:intVsCanc}, a monoid scheme $X$ is integral if and only if it is irreducible, pc and reversible.

\begin{Lem}
\label{lem:VectDecompDual}
Let $X$ be an integral monoid scheme. The functor $F$ of Proposition \ref{prop:VectDecomp} lifts to an exact form functor
\[
(F,\mu):
(
\angleser{\mathcal{O}_X} [\Pic(X)], \mathbf{P}^{\mathcal{L}}, \bm{\Theta}^{\mathcal{L}}) \rightarrow
(\LF^{\nor}(X), P^{\mathcal{L}}, \Theta^{\mathcal{L}}).
\]
\end{Lem}

\begin{proof}
Define the natural isomorphism $\mu: F \circ \mathbf{P}^{\mathcal{L}} \Rightarrow P^{\mathcal{L}} \circ F^{\op}$ so that its components
\[
\mu_{\mathcal{E}_M}: P(\mathcal{E}) \otimes P^{\mathcal{L}}(\mathcal{M}) \rightarrow P(\mathcal{E} \otimes \mathcal{M}) \otimes \mathcal{L},
\qquad
\mathcal{E} \in \LF^{\nor}(X)
\]
are determined by the monoidal data of $P$. It is straightforward to verify that $\mu$ is compatible with $\bm{\Theta}^{\mathcal{L}}$ and $\Theta^{\mathcal{L}}$. We omit the details.
\end{proof}

Define a proto-exact category
\begin{equation}
\label{eq:invtSubcat}
\angleser{\mathcal{O}_X} [\Pic(X)^{P^{\mathcal{L}}}]
=
\bigoplus_{\mathcal{M} \in \Pic(X)^{P^{\mathcal{L}}}} \angleser{\mathcal{O}_X}.
\end{equation}
Since $\Pic(X)^{P^{\mathcal{L}}} \subset \Pic(X)$ is $P^{\mathcal{L}}$-stable, $\angleser{\mathcal{O}_X} [\Pic(X)^{P^{\mathcal{L}}}]$ is a $P^{\mathcal{L}}$-stable proto-exact subcategory of $\angleser{\mathcal{O}_X} [\Pic(X)]$. The induced duality on $\angleser{\mathcal{O}_X} [\Pic(X)^{P^{\mathcal{L}}}]$ is equivalent to that induced by $(P,\Theta)$ on each summand of $\angleser{\mathcal{O}_X}$. In other words, equation \eqref{eq:invtSubcat} defines $\angleser{\mathcal{O}_X} [\Pic(X)^{P^{\mathcal{L}}}]$ as a proto-exact category with duality.

Define a second proto-exact category by
\[
\angleser{\mathcal{O}_X} [\Pic(X)^*\slash P^{\mathcal{L}}]
=
\bigoplus_{\mathcal{M} \in \Pic(X)^* \slash P^{\mathcal{L}}} \angleser{\mathcal{O}_X}.
\]
The choice of a set-theoretic section of the quotient $\Pic(X)^* \rightarrow \Pic(X)^* \slash P^{\mathcal{L}}$ embeds $\angleser{\mathcal{O}_X} [\Pic(X)^*\slash P^{\mathcal{L}}]$ as a proto-exact subcategory of $\LF^{\nor}(X)$ which, however, is not $P^{\mathcal{L}}$-stable. The following result is immediate.

\begin{Prop}
\label{prop:VectDecompDualEnhanced}
Let $X$ be an integral monoid scheme. There is an equivalence of proto-exact categories with duality
\[
\angleser{\mathcal{O}_X} [\Pic(X)]
\simeq
\angleser{\mathcal{O}_X} [\Pic(X)^{P^{\mathcal{L}}}]
\oplus
H \big( \angleser{\mathcal{O}_X} [\Pic(X)^*\slash P^{\mathcal{L}}] \big).
\]
\end{Prop}

\subsection{Grothendieck--Witt theory of monoid schemes}
\label{sec:gwTheoryF1}

Let $\mathcal{L}$ be a line bundle on a reversible pc monoid scheme $X$. Let $\mathcal{GW}(X;\mathcal{L}) = \mathcal{GW}(\LF^{\nor}(X), P^{\mathcal{L}}, \Theta^{\mathcal{L}})$, defined by either via the hermitian $Q$-construction or group completion. Set $GW_i(X;\mathcal{L}) = \pi_i \mathcal{GW}(X;\mathcal{L})$.

\begin{Prop}
Let $X$ be an integral monoid scheme. Then the homotopy type of $\mathcal{GW}(X;\mathcal{L})$ depends on $\mathcal{L}$ only through its class in $\Pic(X) \slash \Pic(X)^2$.
\end{Prop}

\begin{proof}
This follows from Proposition \ref{prop:GWFunctoriality} and Lemma \ref{lem:twistedParity}.
\end{proof}

We have the following analogue of Theorem \ref{thm:KThySchemes} for Grothendieck--Witt theory.

\begin{Thm}
\label{thm:GWThySchemes}
Let $\mathcal{L}$ be a line bundle on an integral monoid scheme $X$. Then there is a natural homotopy equivalence
\[
\mathcal{GW}(X;\mathcal{L})
\simeq
\sideset{}{'}\prod_{\mathcal{M} \in \Pic(X)^{P^{\mathcal{L}}}} \mathcal{GW}(\angleser{\mathcal{O}_X})
\times
\sideset{}{'}\prod_{\mathcal{M} \in \Pic(X)^* \slash P^{\mathcal{L}}} \mathcal{K}(\angleser{\mathcal{O}_X}).
\]
\end{Thm}

\begin{proof}
By Lemma \ref{lem:reflectExactIsoDual}, Proposition \ref{prop:VectDecompDualEnhanced} and the fact that $\mathcal{GW}$ commutes with direct sums of categories, there is a natural homotopy equivalence
\[
\mathcal{GW}(X;\mathcal{L})
\simeq
\mathcal{GW}(\angleser{\mathcal{O}_X} [\Pic(X)^{P^{\mathcal{L}}}]) \times \mathcal{GW}\left(
H(\angleser{\mathcal{O}_X} [\Pic(X)^*\slash P^{\mathcal{L}}]) \right).
\]
The second factor is
\[
\mathcal{GW}\left(
H(\angleser{\mathcal{O}_X} [\Pic(X)^*\slash P^{\mathcal{L}}]) \right)
\simeq
\mathcal{K}(\angleser{\mathcal{O}_X} [\Pic(X)^*\slash P^{\mathcal{L}}])
\simeq
\sideset{}{'}\prod_{\mathcal{M} \in \Pic(X)^* \slash P^{\mathcal{L}}} \mathcal{K}(\angleser{\mathcal{O}_X}),
\]
where the first homotopy equivalence follows from Proposition \ref{prop:GWHyper} and the second from the proof of Theorem \ref{thm:KThySchemes}.

Turning to the factor $\mathcal{GW}(\angleser{\mathcal{O}_X} [\Pic(X)^{P^{\mathcal{L}}}])$, note that since equation \eqref{eq:invtSubcat} respects the duality structures, there is a homotopy equivalence
\[
\mathcal{GW}(\angleser{\mathcal{O}_X} [\Pic(X)^{P^{\mathcal{L}}}])
\simeq
\sideset{}{'}\prod_{\mathcal{M} \in \Pic(X)^{P^{\mathcal{L}}}} \mathcal{GW}(\angleser{\mathcal{O}_X}).
\]
This completes the proof.
\end{proof}

Using Proposition \ref{prop:lineSubCat}, we can combine Theorems \ref{thm:KThyRPC} and \ref{thm:GWThyMonoid} with Theorem \ref{thm:GWThySchemes} to obtain an explicit description of $\mathcal{GW}(X;\mathcal{L})$.

\begin{Cor}
\label{cor:GWIntProp}
Let $X$ be a proper integral monoid scheme. Then there is a homotopy equivalence
\[
\mathcal{GW}^{\oplus}(X;\mathcal{L})
\simeq
\sideset{}{'}\prod_{\mathcal{M} \in \Pic(X)^{P^{\mathcal{L}}}} \mathbb{Z}^2 \times B(\Sigma_{\infty} \times (\mathbb{Z} \slash 2 \wr \Sigma_{\infty}))^+
\times
\sideset{}{'}\prod_{\mathcal{M} \in \Pic(X)^* \slash P^{\mathcal{L}}} \left( \mathbb{Z} \times B \Sigma_{\infty}^+ \right).
\]
\end{Cor}

\begin{proof}
The only additional piece of information needed is Proposition \ref{prop:intPropAut}.
\end{proof}

Without the properness assumption, there is an analogue of Corollary \ref{cor:GWIntProp} is written in terms of $\SPic(\Gamma X)$ and the isometry groups $I(\xi)$. Since we will not use this, we omit its formulation.

Specializing to direct sum (Grothendieck--)Witt groups, we obtain the following results.

\begin{Thm}
\label{thm:GW0Pic}
Let $\mathcal{L}$ be a line bundle on an integral monoid scheme $X$.
\begin{enumerate}[wide,labelwidth=!, labelindent=0pt,label=(\roman*)]
\item There is an isomorphism of abelian groups
\[
GW^{\oplus}_0(X;\mathcal{L})
\simeq
\SPic(\Gamma X) [\Pic(X)^{P^{\mathcal{L}}}] \times \mathbb{Z}[\Pic(X) \slash P^{\mathcal{L}}].
\]

\item There is an isomorphism of abelian groups
\[
W^{\oplus}_0(X;\mathcal{L})
\simeq
\SPic(\Gamma X) [\Pic(X)^{P^{\mathcal{L}}}].
\]
\end{enumerate}
\end{Thm}

\begin{proof}
This follows from Theorem \ref{thm:GWThySchemes}, after using Theorems \ref{thm:KThyRPC} and \ref{thm:GWThyMonoid}. We omit the details.
\end{proof}

\begin{Ex}
Let $\mathbb{A}^1_{\mathbb{F}_1} = \Spec(\mathbb{F}_1[t])$. Since $\mathbb{F}_1[t]$ has no non-trivial idempotents and only the trivial automorphism, the functor $- \otimes_{\mathbb{F}_1} \mathbb{F}_1[t]: \Vect_{\mathbb{F}_1} \rightarrow \LF^{\nor}(\mathbb{A}^1_{\mathbb{F}_1})$ induces an equivalence on maximal groupoids. Moreover, this equivalence respects dualities. It follows that there are homotopy equivalences $\mathcal{K}(\mathbb{A}^1_{\mathbb{F}_1}) \simeq \mathcal{K}(\Vect_{\mathbb{F}_1})$ and $\mathcal{GW}(\mathbb{A}^1_{\mathbb{F}_1}) \simeq \mathcal{GW}(\Vect_{\mathbb{F}_1})$.
\end{Ex}

\begin{Ex}
Let $X = \mathbb{P}_{\mathbb{F}_1}^n$. Fix $d \in \mathbb{Z}$ and set $\mathcal{L} = \mathcal{O}_{\mathbb{P}_{\mathbb{F}_1}^n}(d)$. The involution $P^{\mathcal{L}}$ of $\Pic(X) \simeq \mathbb{Z}$ is $k \mapsto -k + d$. In particular, $\Pic(X)^{P^{\mathcal{L}}}$ is non-empty if and only if $d$ is even, in which case $\Pic(X)^{P^{\mathcal{L}}} = \{\frac{d}{2}\}$. Let $[d]=0$ if $d$ is even and $[d]=1$ otherwise. With this notation, we obtain from Theorem \ref{thm:GW0Pic} an isomorphism
\[
\mathbb{Z}^{[d]} \times \mathbb{Z}[\mathbb{Z}_{\geq d}]
\simeq
GW^{\oplus}_0(\mathbb{P}_{\mathbb{F}_1}^n; d)
,\qquad
b l_{\frac{d}{2}} + \sum_{i \geq d} a_i l_i
\mapsto b [\mathcal{O}_{\mathbb{P}_{\mathbb{F}_1}^n}( \tfrac{d}{2} )] + \sum_{i \geq d} a_i [H(\mathcal{O}_{\mathbb{P}_{\mathbb{F}_1}^n}(i))].
\]
Note that $\mathcal{O}_{\mathbb{P}_{\mathbb{F}_1}^n}( \tfrac{d}{2} )$ admits a unique symmetric form, which is omitted from the notation. In particular, $GW^{\oplus}_0(\mathbb{P}_{\mathbb{F}_1}^n; d)$ is independent of $n$ and, as guaranteed by Lemma \ref{lem:twistedParity}, depends on $d$ only through its parity. Moreover, we have $W^{\oplus}_0(\mathbb{P}_{\mathbb{F}_1}^n; d) \simeq \mathbb{Z}^{[d]}$ with generator $\mathcal{O}_{\mathbb{P}_{\mathbb{F}_1}^n}(\frac{d}{2} )$.
\end{Ex}

\begin{Ex}
Using Theorem \ref{thm:GWThySchemes} and the isomorphism $\Pic(\mathbb{P}^n_{\mathbb{F}_1}) \simeq \mathbb{Z}$, we find
\begin{eqnarray*}
\mathcal{GW}^Q(\mathbb{P}^n_{\mathbb{F}_1})
&\simeq&
\mathcal{GW}^Q(\Vect_{\mathbb{F}_1})
\times
\sideset{}{'}\prod_{k \in \mathbb{Z}_{>0}} \mathcal{K}(\Vect_{\mathbb{F}_1}) \\
& \simeq &
\bigsqcup_{n \in \mathbb{Z}_{\geq 0}}B \Sigma_n \times \mathbb{Z} \times B (\mathbb{Z} \slash 2 \wr \Sigma_{\infty})^+
\times
\sideset{}{'}\prod_{k \in \mathbb{Z}_{>0}} \mathbb{Z} \times B(\Sigma_{\infty})^+.
\end{eqnarray*}
Again, the result is independent of $n$.
\end{Ex}

\subsection{A projective bundle formula}
\label{sec:projBundGW}

We begin with a lemma.

\begin{Lem}
\label{lem:equivarOfPhi}
Let $\pi:\P\cE \to X$ be a projective bundle on an integral monoid scheme. Then the isomorphism
\[
\varphi: \Pic(X) \times \mathbb{Z}\rightarrow \Pic(\mathbb{P}\mathcal{E}),
\qquad
(\mathcal{M},m) \mapsto \pi^\ast (\mathcal{M})\otimes\cO_{\P\cE}(m)
\]
from Theorem \ref{thm:ProjLineBundles} is $\mathbb{Z} \slash 2$-equivariant, where $\mathbb{Z} \slash 2$ acts on $\Pic(X) \times \mathbb{Z}$ by $P^{\mathcal{L}}$ and negation on the first and second factors, respectively, and on $\Pic(\mathbb{P}\mathcal{E})$ by $P^{\pi^*\mathcal{L}}$.
\end{Lem}

\begin{proof}
This is a direct calculation.
\end{proof}

\begin{Thm}
\label{thm:projBundGWSpaces}
Let $\pi:\P\cE \to X$ be a projective bundle on an integral monoid scheme and $\mathcal{L}$ a line bundle on $X$. Then there is a homotopy equivalence
\[
\mathcal{GW}(\mathbb{P}\mathcal{E}; \pi^*\mathcal{L})
\simeq
\mathcal{GW}(X; \mathcal{L}) \times 
\sideset{}{'}\prod_{(\mathcal{M},i) \in (\Pic(X) \times \mathbb{Z}^*) \slash \langle (P^{\mathcal{L}},-1) \rangle} \mathcal{K}(\angleser{\mathcal{O}_X}).
\]
\end{Thm}

\begin{proof}
By Lemma \ref{lem:equivarOfPhi}, the map $\varphi$ induces a bijection
\begin{equation}
\label{eq:fixedBij}
\Pic(\mathbb{P} \mathcal{E})^{P^{\pi^*\mathcal{L}}}
\rightarrow
(\Pic(X) \times \mathbb{Z})^{(P^{\mathcal{L}},-1)}
=
\Pic(X)^{P^{\mathcal{L}}} \times \{0\}
\end{equation}
and a $\mathbb{Z} \slash 2$-equivariant bijection
\begin{eqnarray*}
\Pic(\mathbb{P}\mathcal{E})^*
&\simeq&
\Pic(X)^{P^{\mathcal{L}}} \times \mathbb{Z}^*
\sqcup
\Pic(X)^* \times \mathbb{Z} \\
&=&
\Pic(X)^* \times \{0\}
\sqcup
\Pic(X) \times \mathbb{Z}^*.
\end{eqnarray*}
Together with Theorem \ref{thm:GWThySchemes}, these bijections yield homotopy equivalences
\begin{eqnarray*}
\mathcal{GW}(\mathbb{P}\mathcal{E};\pi^*\mathcal{L})
&=&
\sideset{}{'}\prod_{\mathcal{M} \in \Pic(X)^{P^{\mathcal{L}}}} \mathcal{GW}(\angleser{\mathcal{O}_{\mathbb{P}\mathcal{E}}})
\times\\
&&
\sideset{}{'}\prod_{(\mathcal{M},i) \in (\Pic(X) \times \mathbb{Z}^*) \slash \langle (P^{\mathcal{L}},-1) \rangle} \mathcal{K}(\angleser{\mathcal{O}_{\mathbb{P}\mathcal{E}}})
\times
\sideset{}{'}\prod_{\mathcal{M} \in \Pic(X)^* \slash P^{\mathcal{L}}} \mathcal{K}(\angleser{\mathcal{O}_{\mathbb{P}\mathcal{E}}}) \\
&\simeq&
\mathcal{GW}(X; \mathcal{L}) \times 
\sideset{}{'}\prod_{(\mathcal{M},i) \in (\Pic(X) \times \mathbb{Z}^*) \slash \langle (P^{\mathcal{L}},-1) \rangle} \mathcal{K}(\angleser{\mathcal{O}_{\mathbb{P}\mathcal{E}}}),
\end{eqnarray*}
as claimed.
\end{proof}

Passing to homotopy groups, we obtain $\mathbb{F}_1$-linear analogues of previously known results over fields in which $2$ is invertible \cite{walter2003}, \cite{schlichting2017}, \cite{rohrbach2020}.

\begin{Cor}
\label{cor:projBunGW0}
Let $\pi:\P\cE \to X$ be a projective bundle on an integral monoid scheme and $\mathcal{L}$ a line bundle on $X$.
\begin{enumerate}[wide,labelwidth=!, labelindent=0pt,label=(\roman*)]
\item The map
\[
\varphi: GW_0(X; \mathcal{L}) \times \left( K_0(X) \otimes_{\mathbb{Z}} \mathbb{Z}[\mathbb{Z}^*] \right)^{(\pi^*P^{\mathcal{L}},-1)}
\rightarrow
GW_0(\mathbb{P}\mathcal{E}; \pi^* \mathcal{L})
\]
defined by $\varphi(\mathcal{M}) =\pi^* \mathcal{M}$ and $\varphi(\mathcal{W},m)= H^{\mathcal{L}}(\pi^* \mathcal{W} \otimes_{\mathcal{O}_{\mathbb{P}\mathcal{E}}} \mathcal{O}_{\mathbb{P}\mathcal{E}}(m))$ is an isomorphism of abelian groups.

\item There is an isomorphism of abelian groups
\[
W_0(X; \mathcal{L}) \simeq W_0(\mathbb{P}\mathcal{E}; \pi^* \mathcal{L}).
\] 
\end{enumerate}
\end{Cor}

\begin{proof}
The first statement follows from Theorem \ref{thm:projBundGWSpaces} by taking connected components. Alternatively, we could use Theorems \ref{cor:K0Pic} and \ref{thm:GW0Pic} and Lemma \ref{lem:equivarOfPhi}.

Turning to the second statement, Theorem \ref{thm:GW0Pic} gives
\[
W_0(\mathbb{P}\mathcal{E}; \pi^*\mathcal{L})
\simeq
W_0(\Gamma(\mathbb{P}\mathcal{E},\mathcal{O}_{\mathbb{P}\mathcal{E}}) \proj^{\nor})[\Pic(\mathbb{P}\mathcal{E})^{P^{\pi^*\mathcal{L}}}].
\]
Using Theorem \ref{thm:ProjLineBundles} and the bijection \eqref{eq:fixedBij}, we conclude.
\end{proof}

\begin{Ex}
We have
\[
GW_0(\Spec(\mathbb{F}_1)) \times K_0(\Spec(\mathbb{F}_1)) \otimes_{\mathbb{Z}} \mathbb{Z}[\mathbb{Z}^*]^{\mathbb{Z} \slash 2}
\simeq
\mathbb{Z}^2 \times \mathbb{Z} \otimes_{\mathbb{Z}} \mathbb{Z}[\mathbb{Z}^*]^{\mathbb{Z} \slash 2}
\simeq
\mathbb{Z}^2 \times \mathbb{Z}[\mathbb{Z}_{> 0}]
\]
which is isomorphic to $GW_0(\mathbb{P}_{\mathbb{F}_1}^n) \simeq \mathbb{Z} \times \mathbb{Z}[\mathbb{Z}_{\geq 0}]$, as required by Corollary \ref{cor:projBunGW0}.
\end{Ex}


\appendix
\section{}

\subsection{Direct sums of categories}\label{app:DirectSumCategories}

Let $\{\mathcal{C}_i\}_{i\in I}$ be a family of categories indexed by a set $I$.

\begin{Def}
The direct sum category $\bigoplus_{i\in I} \mathcal{C}_i$ has objects which are finite lists $V_{i_1}\in \mathcal{C}_{i_1}, \dots, V_{i_n}\in \mathcal{C}_{i_n} $ labelled by distinct $i_1, \dots, i_n \in I$. Write $\bigoplus_{j=1}^n V_{i_j}$ for such an object. Morphisms are given by
\[
\Hom_{\bigoplus_{i\in I} \mathcal{C}_i}(\bigoplus_{k=1}^n V_{i_k}, \bigoplus_{l=1}^m W_{j_l})
=
\prod_{\substack{k,l \\ i_k = j_l}} \Hom_{\mathcal{C}_{i_k}}(V_{i_k}, W_{j_l}).
\]
\end{Def}
Many properties and structures of the individual categories $\mathcal{C}_i$ extend in a pointwise fashion to $\bigoplus_{i\in I} \mathcal{C}_i$. For example, if all $\mathcal{C}_i$ are proto-exact (with exact direct sum), then so too is $\bigoplus_{i\in I} \mathcal{C}_i$. 

One can realize $\bigoplus_{i \in I} \mathcal{C}_i$ as a filtered colimit of finite direct sums. Let $\mathscr{P}_{< \infty}(I) $ be the partially ordered set of finite subsets of $I$, ordered by inclusion. Consider the functor $\mathscr{P}_{< \infty}(I) \rightarrow \mathsf{Cat}$ which assigns to a finite subset $S \subset I$ the category $\bigoplus_{s \in S} \mathcal{C}_s$ and to an inclusion $S \hookrightarrow T$ the obvious functor $\bigoplus_{s \in S} \mathcal{C}_s \hookrightarrow \oplus_{t \in T} \mathcal{C}_t$. 

\begin{Prop}\label{prop:DirectSumIsFilteredColimit}
There is an equivalence of categories
\[
\varinjlim_{S\in \mathscr{P}_{< \infty}(I)} \, \bigoplus_{s \in S} \mathcal{C}_s\simeq \bigoplus_{i\in I} \mathcal{C}_i.
\]
\end{Prop}

If $\mathcal{C}=\mathcal{C}_i$ is a constant family of categories with a symmetric bimonoidal structure $(\oplus,\otimes)$ and $(I,\cdot)$ is an abelian group, then we denote the direct sum by
\[
\mathcal{C}[I]:=\bigoplus_{i\in I} \mathcal{C}.
\]
Define a symmetric bimonoidal structure on $\mathcal{C}[I]$ by extending $\oplus$ componentwise and defining $\otimes$ using the convolution product
\[
\bigoplus_a V_{a}\otimes \bigoplus_b W_{b}:=\bigoplus_c(\bigoplus_{\substack{a,b\\ab=c}}V_{a}\otimes W_{b}).
\]

\subsection{Restricted products} \label{app:RestrictedProduct}

Let $\{(Y_i,*_i)\}_{i\in I}$ be family of pointed topological spaces indexed by a set $I$.

\begin{Def}
The restricted product of $\{(Y_i,*_i)\}_{i\in I}$ is
\begin{align*}
\sideset{}{'}\prod_{i\in I} Y_i =\{(y_i)\, |\, y_i\neq *_i \textnormal{ for only finitely many } i\in I\}\subseteq \prod_{i\in I} Y_i
\end{align*}
equipped with the subspace topology.
\end{Def}

The restricted product can be realized as a filtered colimit of finite products as follows.

\begin{Prop}\label{prop:RestrictedProductIsFilteredColimit}
There is a homeomorphism
\[
\varinjlim_{S\in \mathscr{P}_{< \infty}(I)}\, \prod_{s \in S} Y_i\simeq \sideset{}{'}\prod_{i\in I} Y_i.
\]
\end{Prop}

\bibliographystyle{plain}
\bibliography{mybib}
 

\end{document}